\documentclass[12pt, a4paper]{article}
\usepackage{amsmath}
\usepackage{amssymb}
\usepackage{amscd}
\usepackage{amsthm}
\usepackage[totalwidth=17cm, totalheight=25cm]{geometry}
\usepackage{graphicx}

\newtheoremstyle{boldplain}
{9pt}
{9pt}
{\itshape}
{}
{\bfseries}
{.}
{.5em}
{\thmname{#1}\thmnumber{ #2}\thmnote{ (#3)}}%

\newtheoremstyle{bolddefinition}
{9pt}
{9pt}
{}
{}
{\bfseries}
{.}
{.5em}
{\thmname{#1}\thmnumber{ #2}\thmnote{ (#3)}}%

\theoremstyle{boldplain}

\newtheorem{cor}[equation]{Corollary}
\newtheorem{lem}[equation]{Lemma}
\newtheorem{prop}[equation]{Proposition}

\newtheorem{thm}[equation]{Theorem}

\theoremstyle{bolddefinition}

\newtheorem{rem}[equation]{Remark}

\newfont{\bigbf}{cmbx10 scaled\magstep1}
\normalbaselines
\numberwithin{equation}{section}

\def\C{{\mathbb C}}
\def\R{{\mathbb R}}

\def\acts{\curvearrowright}

\def\ora{\overrightarrow}

\def\2pithird{\frac{2\pi}{3}}

\newcommand {\n}[1]{\mathbb{#1}}

\newcommand{\interior}{\operatorname{int}}

\def \eucch {\Delta_{euc}}

\def \R {\mathbb{R}}
\def \polspace {\mathcal{P}_n(X)}

\begin{document}

\title{The generalized triangle inequalities in 
thick Euclidean buildings of rank 2}
\author{Carlos Ramos-Cuevas\footnote{
cramos@mathematik.uni-muenchen.de}}
\date{April 29, 2011}

\maketitle

\begin{abstract}
We describe the set of possible vector valued side lengths of $n$-gons 
in thick Euclidean buildings of rank 2. This set is determined by a 
finite set of homogeneous linear inequalities, which we call the 
{\em generalized triangle inequalities}. 
These inequalities are given in terms of the combinatorics of the 
spherical Coxeter complex associated to the Euclidean building.  
\end{abstract}

\section{Introduction}

Let $X$ be a symmetric space of noncompact type or a 
thick Euclidean building. We are interested in the following 
geometric question:

{\em Which are the possible side lengths of polygons in $X$?}

In this context the appropriate notion of \emph{length} 
of an oriented geodesic segment 
is given by a vector in the Euclidean Weyl chamber
$\eucch$ associated to $X$.
 If $X=G/K$ is a symmetric space, the full invariant of a segment modulo 
the action of $G$ is precisely this vector-valued length since we can 
identify $X\times X/G\cong \eucch$ (cf. \cite{KLM1}). 
For $X$ a Euclidean building the same notion 
of vector-valued length can be defined (cf. \cite{KLM2}).
We denote by $\polspace\subset\eucch^n$ the set of all
possible $\eucch$-valued side lengths of $n$-gons in $X$.

One of the motivations for considering this geometric problem comes
from the following algebraic question.

{\em How are the eigenvalues of two Hermitian matrices related
to the eigenvalues of their sum?} 

This so-called {\em Eigenvalue Problem} goes
back to 1912 when it was already studied by H. Weyl. 
It is closely related to a special case of
the geometric question above, namely, 
for the symmetric space $X=SL(m,\C)/SU(m)$.
We refer to \cite{KLM1} for more information
on the relation between these two questions
and \cite{Fulton} for more history on this problem. 

In \cite{KLM1} and \cite{KLM2} it is shown that the set $\polspace$ 
depends only on the spherical Coxeter complex associated to $X$ 
(i.e.\ on the spherical Weyl chamber $\triangle_{sph}$). 
We will therefore sometimes
refer to $\mathcal{P}_n(\triangle_{sph})$ as the set of side lengths of 
$n$-gons in $X$ a symmetric space or a Euclidean building with 
$\triangle_{sph}$ as spherical Weyl chamber.

For a symmetric space $X=G/K$ the set of possible side lengths has been 
completely determined in \cite{KLM1}:
$\polspace$ is a finite sided convex polyhedral cone and 
it can be described as the solution set of a finite
set of homogeneous linear inequalities in terms of the Schubert calculus 
in the homology of the 
\emph{generalized Grassmannian manifolds} associated to the 
symmetric space $G/K$. It follows, that for a Euclidean building $X'$ with the same
associated spherical Weyl chamber $\triangle_{sph}$ as $X$, the set
$\mathcal{P}_n(X')$ is also a finite sided convex polyhedral cone
determined by the same inequalities as 
$\mathcal{P}_n(X)=\mathcal{P}_n(\triangle_{sph})$.

As already pointed out in \cite{KLM2} for the case of \emph{exotic} 
spherical Coxeter complexes (i.e.\ when it is the Coxeter complex of a 
Euclidean building but it does not occur for a symmetric space) the structure 
of the set $\mathcal{P}_n(\triangle_{sph})$ cannot be described
with this method, since we do not have a Schubert calculus for these
Coxeter complexes. 
Thus, the structure of $\mathcal{P}_n(\triangle_{sph})$
for these Coxeter complexes and even its convexity were unknown.
It is clear that we can restrict our attention to irreducible
Coxeter complexes. 
By a result of Tits \cite{Tits_w}, 
exotic irreducible Coxeter
complexes occur only in rank 2. 
Our main result is the description of $\polspace$
in this case (compare with Theorem~\ref{thm:mainthm}).

\begin{thm}\label{thm:intro}
 For a thick Euclidean building $X$ of rank 2, the space $\polspace$
is a finite sided convex polyhedral cone. The set of inequalities defining
$\polspace$ can be given in terms of the combinatorics of the
spherical Coxeter complex associated to $X$.
\end{thm}

The inequalities given in our main theorem coincide with the so-called
{\em weak triangle inequalities} (cf. \cite[Sec. 3.8]{KLM1}).
Moreover, our arguments also work
(see Theorem~\ref{thm:weakineq}) to prove the  
weak triangle inequalities for buildings of arbitrary rank
(cf. \cite[Thm. 3.34]{KLM1}).
For symmetric spaces,
these inequalities correspond to specially simple intersections of 
Schubert cells in the description of $\polspace$ given in \cite{KLM1}.
Their description depend only in the Weyl group of $X$ and therefore, 
they can be defined for arbitrary Coxeter complexes.

Consider the {\em side length map}
$\sigma: Pol_n(X)=X^n\rightarrow \eucch^n$.
The set $\polspace$ which we are interested in is nothing else than the image of $\sigma$.
We use a direct geometric approach to
describe this image.
Our main idea is to study the singular values 
of $\sigma$ by deforming the sides of a given polygon in $X$.
This strategy was already used for the case of symmetric spaces 
by B.\ Leeb in \cite{Leebnotes}
to give a simple proof of the Thompson Conjecture 
(cf. \cite[Theorem 1.1]{KLM1}). In this paper we adapt this 
{\em variational} method to the case of Euclidean buildings and use it
to describe the space $\mathcal{P}_n(X)$.
This approach is new in the sense that it does not rely
in Schubert calculus at all.

Throughout this paper we state the results, whenever possible,
in such a way that they apply to Euclidean buildings of arbitrary
rank. In particular, 
Sections~\ref{sec:polygons}, \ref{sec:criticalvalues} and
\ref{sec:crosswalls} (except Lemma~\ref{lem:crossablewallsrank2}
and Proposition~\ref{prop:unionofcones})
do not use the assumption on
the rank of the building. 
And when we do use the assumption, we indicate it explicitly 
in the statement of the corresponding result.

The set of inequalities obtained in Theorem~\ref{thm:intro} constitute 
an irredundant system defining the polyhedral cone $\mathcal{P}_n(X)$.
The inequalities given by Schubert calculus in \cite{KLM1} are known
to be irredundant for the cases of type $A_n$ (see \cite{KTW}), however,
these seem to be the only cases. 
A smaller set of inequalities is given
in \cite{BK} by defining a new product in the cohomology of flag varieties.
The irredundancy of this set has been recently shown in \cite{irred}.

After a first version of this paper was written, the author learned about
a recent related paper of Berenstein and Kapovich \cite{BerensteinKapovich},
where the generalized triangle inequalities for rank 2
are also determined by a different approach.

\tableofcontents

\section{Preliminaries}

A very good introduction to the concepts used in this paper is the work
\cite[ch. 2-4]{KleinerLeeb}.
We refer also to \cite{BridsonHaefliger}
for more information on metric spaces with upper curvature bounds and to 
\cite[ch. 2-3]{KLM2} for the different
concepts of {\em length} in Euclidean buildings.

\subsection{CAT(0) spaces}

Recall that a complete geodesic 
metric space $X$ is said to be $CAT(0)$ if the geodesic triangles
in $X$ are not {\em thicker} that the corresponding triangles in the Euclidean
space.

For two points $x,y\in X$ we denote with $xy$ the geodesic segment between them.
The link $\Sigma_x X$ is the completion of the space of directions at $x$ with
the angle metric.
$\ora{xy}\in \Sigma_x X$ denotes the direction of the segment $xy$ at $x$. 

Two complete geodesic lines $\gamma_1,\gamma_2$ are said to be 
{\em parallel} if they have finite Hausdorff distance, 
or equivalently, if the functions
$d(\cdot,\gamma_i)|_{\gamma_{3-i}}$ are constant.
The {\em parallel set}
 $P_\gamma$ is defined as the union of all geodesic lines
parallel to $\gamma$. It is a closed convex set that splits as a metric product
$P_\gamma\cong \R\times Y$, where $Y$ is also a $CAT(0)$ space.

Let $\partial_T X$ denote the Tits boundary of $X$.
For $x\in X$ and $\xi\in \partial_T X$ we write 
$\exp_x(t\xi)$ to denote the point in the geodesic ray $x\xi$ with distance $t$ to $x$.

For a {\em polygon} $p$, or more precisely, 
an {\em $n$-gon} in $X$ we mean the union
of $n$ oriented 
geodesic segments $x_0x_1,\dots,x_{n-1}x_n$ with $x_n=x_0$.
Since geodesic segments in $CAT(0)$ spaces between two given points are unique,
we can also describe $p$ by its {\em vertices}. We write 
$p=(x_0,\dots,x_{n-1})$.
The union $q$ of $n$ oriented 
geodesic segments $x^0x^1,\dots,x^{n-1}x^n$ where
$x^n\neq x^0$ will be called a {\em polygonal path} and we write
$q=(x^0,\dots,x^n)$. We make now the convention that subindices are to
be understood modulo $n$, whereas superindices are not.

\subsection{Coxeter complexes}

A {\em spherical Coxeter complex} is a pair $(S,W)$ consisting of a 
unit sphere $S$ with its usual metric and a finite group $W$ of isometries,
the {\em Weyl group},
generated by reflections on total geodesic spheres of codimension one.
The set of fixed points of reflections in $W$ are called {\em walls}
of $(S,W)$.
A {half-apartment} or {\em root} is a hemisphere bounded by a wall.
A {\em Weyl chamber} in $S$ is a fundamental domain of the action $W\acts S$.
The {\em model Weyl chamber} is defined as $\Delta_{sph}:=S/W$.
We say that two points in $S$ have the same {\em $W$-type} 
(or just {\em type}) if they belong to the same $W$-orbit.

A {\em Euclidean Coxeter complex} is a pair $(E,W_{aff})$
consisting of a Euclidean space $E$ and a group of isometries $W_{aff}$,
the {\em affine Weyl group},
generated by reflections on hyperplanes and such that its rotational part
$W:=rot(W_{aff})$ is finite.
The {\em translation subgroup} $L\subset W_{aff}$ is defined as the kernel 
of the map $rot:W_{aff}\rightarrow W$.
The set of fixed points of reflections in $W_{aff}$ are called {\em walls}
of $(E,W_{aff})$.
A {half-apartment} or {\em root} is a half-space bounded by a wall.
We define the {\em $W_{aff}$-type} of a point in $E$ as 
in the spherical case above.
To $(E,W_{aff})$, we can associate 
the spherical Coxeter complex
$(S, W)$, where $S:=\partial_T E$ is the Tits boundary of $E$.
The {\em Euclidean model Weyl chamber} $\eucch$ 
is the complete Euclidean cone over $\Delta_{sph}$.

The link $\Sigma_x E$ of a point $x\in E$ is naturally a
spherical Coxeter complex with Weyl group $Stab_{W_{aff}}(x)$.

The {\em refined length} of the 
oriented geodesic segment $xy\subset E$ is defined as the image
of $(x,y)$ under the projection 
$E\times E \rightarrow (E\times E)/W_{aff}$.
The {\em $\Delta$-valued length}, or just {\em length},
is the image of the refined length under the natural
{\em forgetful} map $(E\times E)/W_{aff}\rightarrow \eucch$.
We denote with $\sigma$ the {\em length map} assigning to
a segment its $\Delta$-valued length.

We can also define the {\em refined length} of an oriented segment
$xy$ in the spherical Coxeter complex $(S,W)$ analogously as the image
of $(x,y)$ under the projection 
$S\times S \rightarrow (S\times S)/W$.

\subsection{Buildings}

For an introduction to spherical and Euclidean buildings from the point
of view of metric geometry, we refer to \cite{KleinerLeeb}.

Let $X$ be a thick Euclidean building modeled in the Euclidean
Coxeter complex $(E,W_{aff})$.
The concepts of refined length and $\Delta$-valued length
of an oriented geodesic segment $xy\subset X$ can be also defined naturally
by identifying an apartment containing $xy$ with the Coxeter complex
$(E,W_{aff})$.

For a polygon $p=(x_0,\dots,x_{n-1})$ in $X$, we write 
$\sigma(p)=(\sigma(x_0x_1),\dots,\sigma(x_{n-1}x_0))\in\eucch^n$
and call $\sigma:X^n\rightarrow \eucch^n$ the {\em side length map}.
The space $\polspace:=\sigma(X^n)$ is the set of possible
$\Delta$-valued side lengths of $n$-gons in $X$.
We say that a polygon in $X$ is {\em regular}
if all its sides are regular, that is
if their $\Delta$-valued lengths lie in the interior of $\eucch$.
The space of regular polygons is an open dense subset of $X^n$.

We will use following result from \cite{KLM2}
concerning the refined side lengths of polygons in $X$.
We reproduce here its statement for the convenience of the reader.

\begin{thm}[Transfer theorem]\label{thm:transfer}
 Let $X$ and $X'$ be thick Euclidean buildings modeled on
the same Euclidean Coxeter complex $(E,W_{aff})$.
Let $p=(x_0,\dots,x_{n-1})$ be a polygon in $X$ and let $x_0'x_1'$
be a segment in $X'$ with the same refined length as $x_0x_1$.
Then there exists a polygon $p'=(x_0',x_1',\dots,x_{n-1}')$ in $X'$
with the same refined side lengths as of $p$.
\end{thm}

\section{The set of functionals ${\cal L}_n$}\label{sec:functional}

We fix a vertex $o$ of $(E,W_{aff})$ with 
$Stab_{W_{aff}}(o)\cong W$. 
We obtain in this way a natural identification
(modulo the action of $Stab_{W_{aff}}(o)$)
$E\cong \R^{\dim E}$.
By fixing $o$ we get also an embedding $W\hookrightarrow W_{aff}$
and the (coarser) structure $(E,W)$ as Euclidean Coxeter complex.
We will think of the Euclidean Weyl chamber
$\eucch\cong E/W$ as embedded in $E$, such that 
$\eucch$ is a fundamental domain of the action $W\acts E$. 
Hence, the cone point of $\eucch$ corresponds to $o$
and $\partial_T \eucch$ is a Weyl chamber of $(\partial_T E,W)$.

Let $\eta\in \partial_T E$ be a vertex and let $v_\eta:=\exp_o(1\!\cdot\!\eta)$.
We define the following linear functional:
\begin{align*}
l_\eta : \eucch &\rightarrow \n R \\
v &\mapsto \langle ov,ov_\eta \rangle
\end{align*}
where $\langle\cdot,\cdot\rangle$ denotes the standard scalar
product on $\R^{\dim E}$.
We denote with ${\cal L}_n$ the finite set of functionals on $\eucch^n$ 
of the form 
$L(v_0,\dots,v_{n-1})=
l_{\eta_0}(v_0)+\dots+l_{\eta_{n-1}}(v_{n-1})$ 
where all the $\eta_i$ have the same $W$-type.
We write $L=(l_{\eta_0},\dots,l_{\eta_{n-1}})$ for such a functional.

Let $H_L$ denote the hyperplane $L^{-1}(0)\cap \eucch^n$ for
$L\in{\cal L}_n$.
We call $H_L$ a {\em wall} in $\eucch^n$.
The set of walls $H_L$ divide $\eucch^n$ in finitely many
convex polyhedral cones. We denote with ${\cal C}_n$ the family of 
the interiors of these 
cones, i.e.\ ${\cal C}_n$ is the set of the connected components
of $\interior(\eucch^n) \setminus \bigcup_{L\in {\cal L}_n}H_L$.

\section{Polygons}\label{sec:polygons}

\subsection{Holonomy map}\label{sec:holonomy}

Let $p=(x_0,\dots, x_{n-1})$ be an $n$-gon in $X$. We say that a 
$n$-tuple ${\cal F} = (F_0, \dots, F_{n-1})$ of apartments in $X$ {\em supports
the polygon} $p$ if $x_i x_{i+1}\subset F_i$ and 
the convex set $F_i\cap F_{i+1}$
is top dimensional and contains $x_{i+1}$ in its interior.

\begin{rem}
 If $p$ is a regular polygon then there always exists an $n$-tuple
${\cal F}$ supporting $p$. ${\cal F}$ can be constructed as follows:
Let $A\in\Sigma_{x_0}X$ be an apartment containing $\ora{x_0x_1}$ and 
$\ora{x_0x_{n-1}}$ and take $\xi\in A$ antipodal to $\ora{x_0x_1}$. 
Extend the segment $x_0x_1$ a little further than $x_0$ in direction of $\xi$ to a segment $x_0'x_1$.
Inductively for $i=0,\dots,{n-2}$ choose $F_i\in X$ to be an 
apartment containing $x_i'x_{i+1}$ and an initial part of $x_{i+1}x_{i+2}$
and extend $x_{i+1}x_{i+2}$ in $F_i$ a little further than $x_{i+1}$ 
to a segment $x_{i+1}'x_{i+2}$.
Finally choose $F_{n-1}$ to contain $x_{n-1}'x_0$, 
an initial part of $x_{0}x_1$ and $x_{0}x_0'$. 
This last step is possible because of our first choice of $x_0'$.
The polyhedron $F_i\cap F_{i+1}$ contains a regular segment with $x_{i+1}$
in its interior. In particular, $F_i\cap F_{i+1}$ is top dimensional.
\end{rem}

Let now $p$ be a polygon and ${\cal F}$ an $n$-tuple supporting it. 
Notice that since the convex set $F_i\cap F_{i+1}$ is 
top dimensional, there is a unique isomorphism of Coxeter complexes
$\phi_i:F_i\rightarrow F_{i+1}$ fixing $F_i\cap F_{i+1}$ pointwise.
We also write $\phi_i$ for the induced isomorphism at the boundary:
$$
\phi_i : S_i:=\partial_T F_i \rightarrow S_{i+1}:=\partial_T F_{i+1}.
$$

We obtain an associated \emph{holonomy map} 
$\phi_p:S_i \rightarrow S_i$ defined as the 
composition $\phi_p:=\phi_{i+n-1}\circ\dots\circ\phi_{i+1}\circ\phi_i$.
We introduce also the following notation:
$$
\phi_i^k:=\phi_{i+k-1}\circ\dots\circ\phi_i:F_i\rightarrow F_{i+k}
$$
and the analogous for the Tits boundaries.
Recall our convention that subindices are taken modulo $n$ and superindices are not.

Notice that the holonomy map $\phi_p:S_i\rightarrow S_i$ 
is an element of the Weyl group $W$. 
In particular the set of fixed points of $\phi_p$
is a singular sphere in $(S_i,W)$. 
We point out that the holonomy map (and therefore also its fixed points set) 
depends on the choice of the $n$-tuple ${\cal F}$ supporting $p$.
We will make use of this flexibility later.

\subsection{Opening a polygon in an apartment}\label{sec:openpolygon}

\begin{figure}[h]
\begin{center}
\includegraphics[scale=0.5]{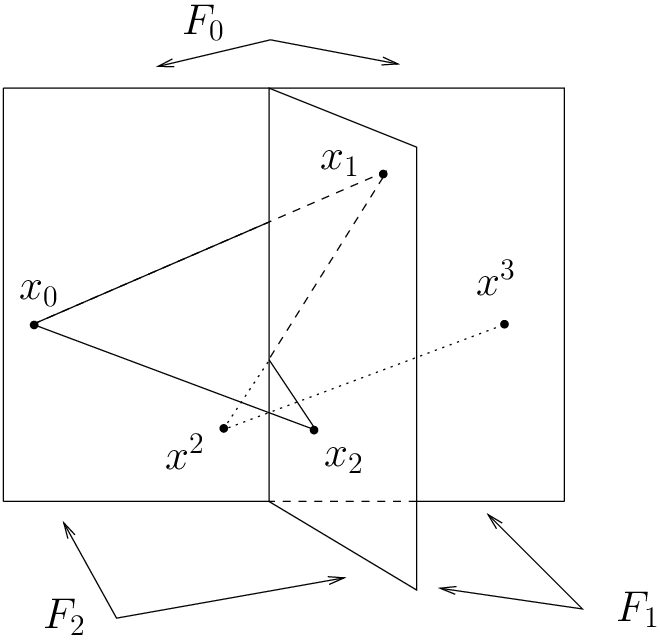}
\caption{Opening a triangle in the apartment $F_0$}\label{fig:openpolygon}
\end{center}
\end{figure}

Let $p=(x_0,\dots, x_{n-1})$ be an $n$-gon in $X$ and let ${\cal F}$
be an $n$-tuple supporting it.
We want to construct a special polygonal path $\bar p = (x^0,\dots,x^n)$
contained in the apartment
$F_0$ which has the same refined side lengths as $p$ and 
{\em looks like} $p$ at the vertices, that is,
such that the segments $\ora{x_ix_{i-1}}\ora{x_ix_{i+1}}$
and $\ora{x^ix^{i-1}}\ora{x^ix^{i+1}}$ have the same refined side lengths
(see Fig.~\ref{fig:openpolygon}).
For this we just take $x^0:=x_0$ and for $i>0$
we define $x^i:=(\phi_0^i)^{-1}(x_i)=(\phi_0^{i-1})^{-1}(x_i)\in F_0$.
Recall that $\phi_i|_{F_i\cap F_{i+1}}=Id$.
The polygonal path $\bar p$ so defined has the required properties since
$x^ix^{i+1}=(\phi_0^i)^{-1}(x_ix_{i+1})$ 
and in $\Sigma_{x^i}F_0$ we have
$\ora{x^ix^{i-1}}=(\phi_0^i)^{-1}(\ora{x_ix_{i-1}})$ and
 $\ora{x^ix^{i+1}}=(\phi_0^i)^{-1}(\ora{x_ix_{i+1}})$.

We remark that in general
$x^n\neq x^0$ and $(x^0,\dots,x^n)$
is a polygonal path
hence the expression ``opening a polygon''.

\subsection{Folding a polygon into an apartment}

This construction was first considered in \cite[Sec.\ 6.1]{KLM3}.

\begin{figure}[h]
\begin{center}
 \includegraphics[scale=0.8]{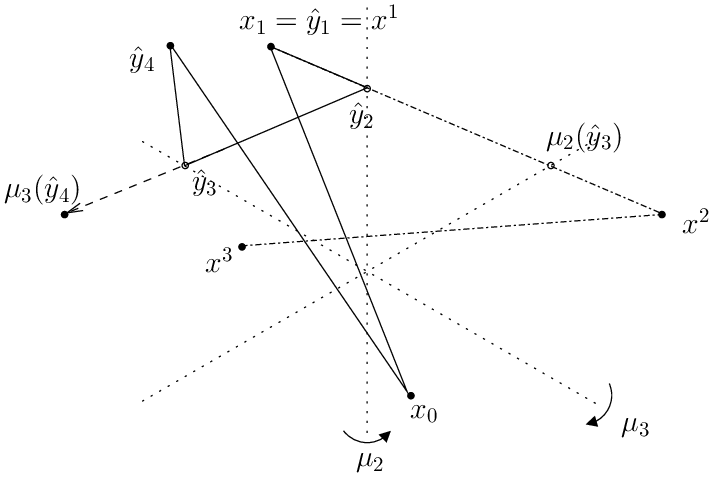}
\caption{Folding and opening a triangle}\label{fig:foldpolygon}
\end{center}
\end{figure}

For simplicity on the notation, suppose $p=(x_0,x_1,x_2)$
is a triangle in $X$. 
There is a partition 
$y_1=x_1,y_1,\dots,y_k=x_2$ 
of the segment $x_1x_2$
such that the triangles
$(x_0,y_i,y_{i+1})$ for $i=1,\dots,k-1$ are contained in an apartment $A_i$,
in particular, they are flat.
We define points $\hat y_i$ in the apartment $A_1$ inductively as follows:
for $i=1$ set $\hat y_1 = y_1 = x_1$ and suppose we have already defined
$\hat y_i$.
Let $\beta_i:A_i\rightarrow A_1$ be an isomorphism of Euclidean Coxeter
complexes, such that $\beta_i(x_0y_i)=x_0\hat y_i$.
We define $\hat y_{i+1}:=\beta_i(y_{i+1})$.
We say that the polygon $\hat p=(x_0,\hat y_1,\dots,\hat y_k)$ 
is the result of {\em folding} the triangle $p$ into $A_1$ (see Fig.~\ref{fig:foldpolygon}). 
We say that the points $\hat y_i$ for $i=2,\dots,k-1$ are
the {\em break points} of the folded polygon $\hat p$.
Notice that the segments $x_0x_1$ and $x_0x_2$ have the same refined
side lengths as the segments $x_0\hat y_1$ and $x_0\hat y_k$ respectively.
Write $\hat y_0=y_0 = x_0$ and define 
$\zeta_i:=\ora{y_iy_{i-1}}$ and $\xi_i:=\ora{y_iy_{i+1}}$, analogously
$\hat\zeta_i:=\ora{\hat y_i\hat y_{i-1}}$ and 
$\hat\xi_i:=\ora{\hat y_i\hat y_{i+1}}$.

A {\em billiard triangle} is a polygon $\hat p=(x_0,\hat y_1,\dots,\hat y_k)$
in an apartment $A$ such that for $i=2,\dots,k-1$ the directions
$\hat\zeta_i$ and $\hat\xi_i$ are antipodal in 
the spherical Coxeter complex $(\Sigma_{\hat y_i}A_1, Stab_{W_{aff}}(\hat y_i))$ 
modulo the action of the Weyl group
$Stab_{W_{aff}}(\hat y_i)$. Clearly, a folded triangle is a billiard triangle.
Conversely, 
the next condition is necessary and sufficient for a billiard triangle to
be a folded triangle.

{\em
For $i=2,\dots,k-1$ there is a triangle $(\zeta_i',\xi_i',\tau_i')$
in the spherical building $\Sigma_{\hat y_i}X$ such that 
$d(\zeta_i',\xi_i')=\pi$ and the refined lengths of
$\zeta_i'\tau_i'$ and $\xi_i'\tau_i'$ are the same as of
$\hat\zeta_i\ora{\hat y_ix_0}$ and
$\hat\xi_i\ora{\hat y_ix_0}$ respectively.
}

We investigate now the relation between the constructions of opening and
folding a polygon. Let $p=(x_0,x_1,x_2)$
be a triangle in $X$ and let ${\cal F}$ be a triple supporting $p$.
Observe that we can choose $A_1=F_0$, $A_{k-1}=F_2$ and $\beta_{k-1}=\phi_2$.
Let $\hat p=(\hat y_0,\hat y_1,\dots,\hat y_k)$ be the folded triangle
and let $\bar p=(x^0,x^2,x^3)$ be the opened triangle
with $\hat y_0 = x_0 = x^0$ and $\hat y_1 = x_1 = x^1$.

Since 
$\hat p=(\hat y_0,\hat y_1,\dots,\hat y_k)$ is a billiard triangle,
there are isometries $\mu_i$ of $F_0$ in $Stab_{W_{aff}}(\hat y_i)$
for $i=2,\dots,k-1$ such that $\ora{\hat y_i\mu_i(\hat y_{i+1})}$ and 
$\hat\zeta_i$ are antipodal in $\Sigma_{\hat y_i} F_0$ (Fig.~\ref{fig:foldpolygon}).
We call the $\mu_i$
the {\em straightening} isometries. It holds
$$
\mu_2\circ\dots\circ\mu_{k-2}\circ\mu_{k-1}(x_0) = x^3\;
\text{ and }\;
\mu_2\circ\dots\circ\mu_{k-2}\circ\mu_{k-1}(\hat y_k) = x^2.
$$

Observe that if $p$ is regular, then the $\mu_i$'s are unique.

The constructions for $n$-gons ($n>3$) are analogous.

\section{Critical values of the side length map $\sigma$}\label{sec:criticalvalues}

For a {\em regular} value of the side length map $\sigma$ we mean
a value $s\in \polspace$ for which there is a polygon $p$ with $\sigma(p)=s$
and such that $\sigma$ is an open map at $p$.
First we give a sufficient condition in terms of the holonomy map
for $\sigma(p)$ being a {\em regular} value of $\sigma$.

\begin{prop}\label{prop:regularvalue}
 Let $p$ be an $n$-gon in $X$ and ${\cal F}$ an $n$-tuple supporting $p$.
Suppose that the holonomy map $\phi_p:S_i\rightarrow S_i$ has no fixed points, then 
the space $\polspace$ is a neighborhood of $\sigma(p)$ in $\eucch^n$.
\end{prop}
\proof
Choose $\epsilon>0$ so that 
$B_{x_i}(n\epsilon)\cap F_{i-1} = B_{x_i}(n\epsilon)\cap F_{i} \subset F_{i-1}\cap F_{i}$
for all $i$. 
Recall that for $x\in X$ and $\xi\in \partial_T X$,  
$\exp_x(t\xi)$ denotes the point in the geodesic ray $x\xi$ with distance $t$ to $x$.
Notice that for $0\leq t < n\epsilon$ and $\xi\in S_i$,
$\exp_{x_i}(t\xi)\in F_{i-1}\cap F_{i}$. The segments 
$x_{i}\exp_{x_{i}}(t\xi)$ and $x_{i+1}\exp_{x_{i+1}}(t\phi_i(\xi))$ are two parallel segments
of the same length in the apartment $F_i$.

\begin{figure}[h]
\begin{center}
 \includegraphics[scale=0.7]{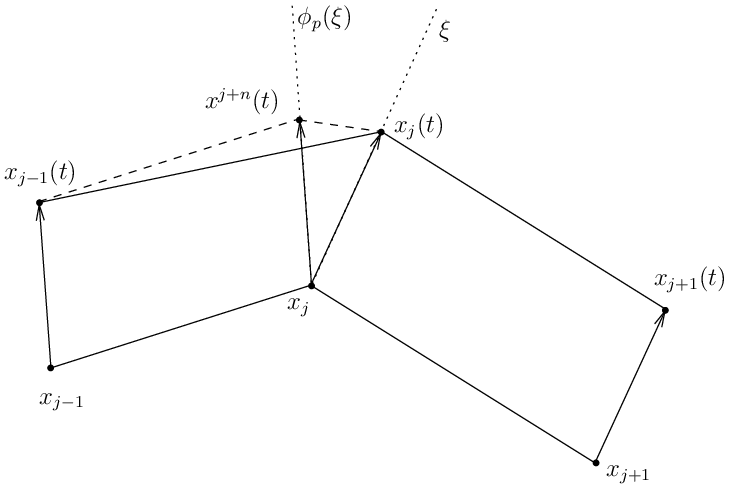}
\caption{Variation of the side $x_{j-1}x_j$}\label{fig:variationofside}
\end{center}
\end{figure}

Fix a $0\leq j \leq n-1$ and $\xi \in S_j$.
We want to {\em move} the polygon $p$ {\em along} the direction $\xi\in S_j$ 
in the following sense.
Let $t<\epsilon$ and set 
$x^{j+k}(t):=\exp_{x_{j+k}}(t\,\phi_j^k(\xi))$ for $k\geq 0$.
Consider the polygon 
$p(\xi,t):=(x_j(t),\dots,x_{j+n-1}(t)):=(x^j(t),\dots,x^{j+n-1}(t))$.
Notice that for $i\neq j$ the segment $x_{i-1}(t)x_i(t)$ is just a
translation in the apartment $F_i$ of the segment $x_{i-1}x_i$,
in particular, they have the same $\Delta$-valued length.
But since the holonomy map has no fixed points, then
$\phi_j^n(\xi)=\phi_p(\xi)\neq \xi$ and we get (see Fig.~\ref{fig:variationofside})
$$
\sigma(x_{j-1}(t)x_j(t))=\sigma( x_{j-1}(t)x^{j+n}(t) + x^{j+n}(t)x_j(t)).
$$

If we think of $x^{j+n}(t)x_j(t)$ as vectors in the Euclidean space $F_{j-1}$,
then the set $\{ x^{j+n}(t)x_j(t) \; | \; \xi\in S_j,\; 0\leq t<n\epsilon\} $
is a neighborhood of the origin, because $\phi_p$ has no fixed points.
Hence, the set $\{\sigma(x_{j-1}(t)x_j(t)) \; | \; \xi\in S_j,\; 0\leq t<n\epsilon\}$ 
is a neighborhood
of $\sigma(x_{j-1}x_j)$ in $\eucch$. 
This means that we can deform the $\Delta$-valued length of every side
of $p$ independently, thus 
$\polspace$ is a neighborhood of $\sigma(p)$ in $\eucch^n$.
\qed

The next proposition says that for a building with only one vertex 
the critical values of $\sigma$ must lie in
the walls $H_L$.

\begin{prop}\label{prop:criticalvalue}
 Let $p$ be an $n$-gon in a thick Euclidean building $X$ which has only one
vertex. Let ${\cal F}$ be an $n$-tuple supporting $p$.
Suppose that the holonomy map $\phi_p:S_i\rightarrow S_i$ 
fixes a vertex of $(S_i,W)$. 
Then there exists a functional $L\in {\cal L}_n$, such that
$L(\sigma(p))=0$.
\end{prop}
\proof
Let $o\in X$ be the only vertex. We have that $W_{aff}\cong W$ and
$W_{aff}\acts F$ fixes $o\in F$ for any apartment $F\subset X$.
Let $\eta\in S_0$ be a vertex fixed by 
$\phi_p:S_0\rightarrow S_0$.
Now open the polygon $p=(x_0,\dots,x_n)$ in the apartment $F_0$
to the polygonal path $\bar p=(x^0,\dots,x^n)$. 
Recall that $x_0=\phi_0^n(x^n)$ and $\phi_0^n=\phi_p\in W_{aff}$.
Thus, $\phi_p(o)=o$.

Let $v:=\exp_o(1\!\cdot\!\eta)$
and let $\eta_i\in \partial_T E$ for $i=0,\dots,n-1$
be vertices of the same $W$-type
as $\eta$, such that for $v_i:=\exp_o(1\!\cdot\!\eta_i)$ holds
$l_{\eta_i}(\sigma(x_{i}x_{i+1}))=\langle x^ix^{i+1}, ov\rangle$.
Set $L=(l_{\eta_1},\dots,l_{\eta_n})$, then
$$
L(\sigma(p))=\int_{\bar p}\langle o\cdot,ov\rangle
=\langle ox^n,ov\rangle-\langle ox_0,ov\rangle
=\langle ox^n,ov\rangle-\langle \phi_p(ox^n),\phi_p(ov)\rangle=0.
$$
\qed

\begin{rem}\label{rem:moreonevertex}
The Figure~\ref{fig:moreonevertex} shows an example of a folded triangle in
a thick Euclidean building 
with more than one vertex, 
which comes from a genuine triangle $p$ (cf.~\cite[Sec. 6.1]{KLM3}).
The holonomy map $\phi_p:S_0\rightarrow S_0$ fixes the direction $\eta$ but
$\phi_p:F_0\rightarrow F_0$ has no fixed points. It follows that for
the associated functional as in the proof of Proposition~\ref{prop:criticalvalue}
holds $L(\sigma(p))\neq 0$. Moreover, if the vertices of $p$ are in
a sufficiently general position, one can achieve that $\sigma(p)$ does not lie in any of
the walls $H_L$. But since the vertices lie in the interior of Weyl alcoves any small
variation of $p$ will leave $L\circ\sigma$ constant. 
Hence, $\sigma$ cannot be open at $p$. 
Nevertheless, we have Corollary~\ref{cor:imageisopen} below.
\begin{figure}[h]
\begin{center}
 \includegraphics[scale=0.7]{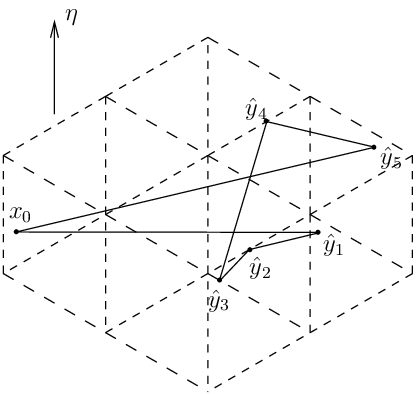}
\caption{A triangle with a direction fixed by the holonomy map but with $L(\sigma(p))\neq 0$}
\label{fig:moreonevertex}
\end{center}
\end{figure}
\end{rem}

We use the result in \cite{KLM2} that $\polspace$ depends only on
the spherical Coxeter complex to transfer the result above to arbitrary
buildings.

\begin{cor}\label{cor:imageisopen}
 Let $s\in\polspace\cap \interior\eucch^n$ and suppose that $L(s)\neq 0$
for all functionals $L\in{\cal L}_n$. Then $\polspace$ is a neighborhood 
of $s$ in $\eucch^n$.
\end{cor}
\proof
By \cite{KLM2} we may assume that $X$ has only one vertex. Let $p$ be a regular
polygon with $\sigma(p)=s$ and let ${\cal F}$ be an $n$-tuple supporting $p$.
By Proposition~\ref{prop:criticalvalue} the holonomy map has no fixed points.
The result now follows from Proposition~\ref{prop:regularvalue}.
\qed

\begin{lem}\label{lem:imageisclosed}
 Let ${p_k}$ be a sequence of regular $n$-gons in $X$ such that 
$\sigma(p_k)\rightarrow s$ in $\eucch^n$, then there exists an $n$-gon
$p$ in $X$ such that $\sigma(p)=s$.
\end{lem}
\proof
We assume again that $X$ has only one vertex $o\in X$. 
Let $p_k=(x_0^k,\dots,x_{n-1}^k)$
and let ${\cal F}_k=(F_0^k,\dots,F_{n-1}^k)$ be $n$-tuples supporting $p_k$.
After transferring the polygons $p_k$ 
(cf. Theorem~\ref{thm:transfer})
we may assume that the sides $x_0^kx_1^k$ lie in the same apartment $F=F_0^k$.
and that $x_0^k$
lie in the same Euclidean Weyl chamber $\eucch\subset F$. 
We open now the polygons $p_k$ in the apartment $F$ to polygonal paths
$\bar p_k=(x^{k,0},\dots,x^{k,n})$.

If $x^{k,0}\rightarrow\infty$ in $F$, then 
after taking a subsequence the segments $ox^{k,0}$ converge to a geodesic ray $\rho$.
Let $\tau\subset \eucch$ be the smallest face of $\eucch$ containing $\tau$.
Then for $k$ big enough, folding the polygon $p_k$ into $F$ 
will have break points only on open faces of $F$ containing $\tau$
in their closure. This implies that 
$p_k$ is contained in the 
parallel set $P_{f_\tau}$ of the flat $f_\tau\subset F$ spanned by $\tau$.
For instance, if $\tau=\eucch$, then $p_k$ is contained in the
apartment $F$ for $k$ big enough.
Thus, after translating the polygons $p_k$ in the parallel set $P_{f_\tau}$ 
we may assume that $x^{k,0}$ stay in a bounded region.

Now we can take a subsequence and assume that the polygonal paths $\bar p_k$ 
converge to a polygonal path $\bar p=(x^0,\dots,x^n)$ with 
$\Delta$-valued side lengths $s$.
We want now to lift this polygonal path near the polygons $p_k$. 
Write $\varphi^{k,i}:F \rightarrow F_i^k$ 
for the maps $\phi_0^i:F_0^k\rightarrow F_i^k$ defined in Section~\ref{sec:holonomy}.
For $i=1,\dots,n$ the convex set 
$(\varphi^{k,i})^{-1}(F_{i-1}^k\cap F_i^k)=(\varphi^{k,i-1})^{-1}(F_{i-1}^k\cap F_i^k)$
is a union of Euclidean Weyl chamber with $x^{k,i}$ in its interior.
Hence, for $k$ big enough we have 
$x^i\in(\varphi^{k,i})^{-1}(F_{i-1}^k\cap F_i^k)=(\varphi^{k,i-1})^{-1}(F_{i-1}^k\cap F_i^k)$
and we can define the points 
$z^{k,i}:=\varphi^{k,i}(x^i)=\varphi^{k,i-1}(x^i)\in F_{i-1}^k\cap F_{i}^k$
for $i=1,\dots,n$ and $z^{k,0}:=x^0$.
Then $q_k:=(z^{k,0},\dots,z^{k,n})$ is a polygonal path with the same side lengths
as $\bar p$, i.e.\ $\sigma(q_k)=s$. However $q_k$ may still not be a closed polygon. 

Notice that 
$d(x^{k,0},z^{k,n})=d(\varphi^{k,n}(x^{k,n}),\varphi^{k,n}(x^{n}))=
d(x^{k,n},x^n)\rightarrow 0$
and $d(x^{k,0},x^0)\rightarrow 0$, thus 
$d(z^{k,n},z^{k,0})=d(z^{k,n},x^0)\rightarrow 0$.
On the other hand,
observe that $x^{k,n}$ and $x^{k,0}=\varphi^{k,n}(x^{k,n})$
have the same $W_{aff}$-type and therefore
also their limits $x^0=z^{k,0}$ and $x^n=(\varphi)^{-1}(z^{k,n})$.
Hence, $z^{k,0}$ and $z^{k,n}$ have the same type. 
But $W_{aff}$ is finite, so $d(z^{k,0},z^{k,n})$
can only take finitely many values. It follows that for $k$ big enough
$z^{k,0}=z^{k,n}$ and $q_k$ is a closed polygon in $X$ 
with $\Delta$-valued side lengths $s$.
\qed

\begin{cor}\label{cor:opencones}
For any open cone $C\in{\cal C}_n$ the intersection $\polspace\cap C$ is 
empty or $C$.
Moreover, if $C\subset \polspace$, then $\bar C\subset \polspace$.
\end{cor}
\proof
The intersection $\polspace\cap C$ is open by Corollary~\ref{cor:imageisopen}
and closed by Lemma~\ref{lem:imageisclosed}.
The second assertion also follows from Lemma~\ref{lem:imageisclosed}.
\qed

\section{The generalized triangle inequalities}

\subsection{Crossing the walls $H_L$}\label{sec:crosswalls}

Suppose $p$ is a polygon in $X$ with $\sigma(p)=s\in H_L$ for some functional
$L\in{\cal L}_n$. 
Considering Corollary~\ref{cor:opencones}
the natural question is if there is a cone $C\in {\cal C}_n$ such that
$s\in \bar C \subset \polspace$. We would also like to describe all cones
in ${\cal C}_n$ with this property. With this in mind we investigate in this section 
following question. When can we find polygons $p'$ with 
$\Delta$-valued side lengths near $s$ and such that $L\circ\sigma(p)>0$ 
(or $<0$)? For this we might try to study the side lengths of small perturbations
of $p$. However since a Euclidean building has dimension equal to his rank,
we do not have much flexibility to perturbate the polygon
(compare with Remark~\ref{rem:moreonevertex}). 
Thus we must be
more compliant with the variations of $p$ that we want to admit.
Therefore we will often have to translate the polygon to other place in $X$
where we can perform the perturbations. 

Let $L=(l_{\eta_0},\dots,l_{\eta_{n-1}})$ be a functional in ${\cal L}_n$.
For the rest of this section 
$p=(x_0,\dots,x_{n-1})$ will be always a regular $n$-gon
such that $\sigma(p)\in H_L$, that is, $L(\sigma(p))=0$.

Let ${\cal F}$ be an $n$-tuple of apartments
supporting $p$.
Let $\xi_i\in S_i$ be a vertex such that if 
$v_i:=\exp_{x_i}(1\!\cdot\!\xi_i)$, then
$l_{\eta_i}(\sigma(x_ix_{i+1}))=\langle x_ix_{i+1}, x_iv_i\rangle$.
Observe that $v_i$ is of the same $W$-type as 
$\eta_0,\dots,\eta_{n-1}$. 
We will therefore sometimes write  
$l_{\eta_i}\circ\sigma=\langle \cdot, x_iv_i\rangle$.

\begin{lem}\label{lem:wallbyaccident}
If in the notation above $\ora{x_i\xi_i}\neq \ora{x_i\xi_{i-1}}$ for some $i$, 
then for any neighborhood
$U$ of $\sigma(p)$ in $\eucch^n$ there exist $n$-gons $p_1,p_2$ in $X$
with $\sigma(p_1),\sigma(p_2)\in U$ and $L\circ\sigma(p_1)>0>L\circ\sigma(p_2)$.
\end{lem}
\proof
The proof is similar to the one of Proposition~\ref{prop:regularvalue}.
Notice that for small $\epsilon>0$ holds
$x_i':=\exp_{x_i}(\epsilon\xi_{i-1}),
x_i'':=\exp_{x_i}(\epsilon\xi_{i}) \in F_{i-1}\cap F_i$
and by the hypothesis $x_i'\neq x_i''$. 
Let $\theta:=\angle_{x_i}(\xi_{i-1},\xi_i)$.
Consider first the polygon $p_1:=(x_0,\dots,x_i',\dots,x_{n-1})$, then
\begin{eqnarray*}
L(\sigma(p_1)) & = & l_{\eta_0}(\sigma(x_0x_1))+\dots+
\langle x_{i-1}x_i+x_ix_i', x_{i-1}v_{i-1} \rangle + 
\langle x_ix_{i+1}-x_ix_i', x_iv_i \rangle +\\
&&\dots + l_{\eta_{n-1}}(\sigma(x_{n-1}x_0)) \\
& = & l_{\eta_0}(\sigma(x_0x_1))+\dots+
(l_{\eta_{i-1}}(\sigma(x_{i-1}x_i)) + \epsilon) + 
(l_{\eta_{i}}(\sigma(x_{i}x_{i+1})) - \epsilon\cos\theta) +\\
&&\dots + l_{\eta_{n-1}}(\sigma(x_{n-1}x_0))\\
& = & L(\sigma(p)) + \epsilon(1-\cos\theta) 
\;\;>\;\; L(\sigma(p)) = 0.
\end{eqnarray*}

Analogously for the polygon $p_2:=(x_0,\dots,x_i'',\dots,x_{n-1})$
we have $L(\sigma(p_2))=L(\sigma(p)) + \epsilon(\cos\theta-1)<L(\sigma(p))=0$.
\qed

Suppose now that $\ora{x_i\xi_i} = \ora{x_i\xi_{i-1}}$ for all $i$. 
In particular, $\phi_{i-1}(\xi_{i-1})=\xi_i$ and the holonomy
map $\phi_p:S_i\rightarrow S_i$ has the fixed point $\xi_i$.
Let $\xi_i'\in S_i$ be the antipodal point to $\xi_i\in S_i$.
If all $\xi_i$ and $\xi_i'$ coincide, 
then the polygon $p$ is contained
in a parallel set, namely the set $P_{\xi_0,\xi_0'}$ of all lines 
connecting $\xi_0$ with $\xi_0'$.

\begin{lem}\label{lem:localparallelset}
Suppose $p$ is not contained in any parallel set $P_{\xi,\xi'}$, 
where $\xi,\xi'\in \partial_T X$ are antipodal points
such that 
$\ora{x_i\xi_i}=\ora{x_i\xi}$ for all $i$.
Then for any neighborhood
$U$ of $\sigma(p)$ in $\eucch^n$ there exist $n$-gons $p_1,p_2$ in $X$
with $\sigma(p_1),\sigma(p_2)\in U$ and $L\circ\sigma(p_1)>0>L\circ\sigma(p_2)$.
\end{lem}
\proof
Let $P=(\nu_0,\dots,\nu_{n-1})$ be an $n$-tuple of geodesic segments
$\nu_i:[s^-,s^+]\rightarrow X$ with 
$s^-<0<s^+$, $\nu_i(0)=x_i$,
$\dot \nu_i(0) = \ora{x_i\xi_i}$
and such that the convex hull $CH(\nu_i,\nu_{i+1})$ is a (2-dimensional)
flat quadrilateral. 
Notice that $CH(\nu_i,\nu_{i+1})$ is actually a parallelogram because
$\ora{x_i\xi_i} = \ora{x_i\xi_{i-1}}$.
Such a $P$ exists, just take small parts of the lines 
through $x_i$ connecting $\xi_i$ with $\xi_i'$.
Suppose now that $P$ is maximal with these properties,
i.e.\ the segments $\nu_i$ cannot be extended.
If $|s^\pm|=\infty$, then the $\nu_i$ are parallel geodesic lines and
$p\subset P_{\nu_0}$. 
Hence at least one
of $s^+$ or $-s^-$ must be $<\infty$. 
Suppose $s=s^+<\infty$ (the other
case is analogous).

Now we want to displace $p$ along $\nu_i$ to the region, 
where it does not look locally like a parallel set anymore: set 
$p'=(x_0',\dots,x_{n-1}')=(\nu_0(s),\dots,\nu_{n-1}(s))$.
Then $p'$ is an $n$-gon with $\sigma(p')=\sigma(p)$.
Choose apartments $A_i$ containing the convex sets $CH(\nu_{i},\nu_{i+1})$.
Let
$\zeta_i:=-\dot{\nu_i}(s)\in\Sigma_{x_i'}(A_{i-1}\cap A_{i})$
and let
$\alpha_i\in\Sigma_{x_i'}A_{i-1}$, 
$\beta_i\in\Sigma_{x_i'}A_{i}$ be the antipodes of $\zeta_i$
in $\Sigma_{x_i'}A_{i-1}$ and $\Sigma_{x_i'}A_{i}$ respectively.

\begin{figure}[h]
\begin{center}
 \includegraphics[scale=0.5]{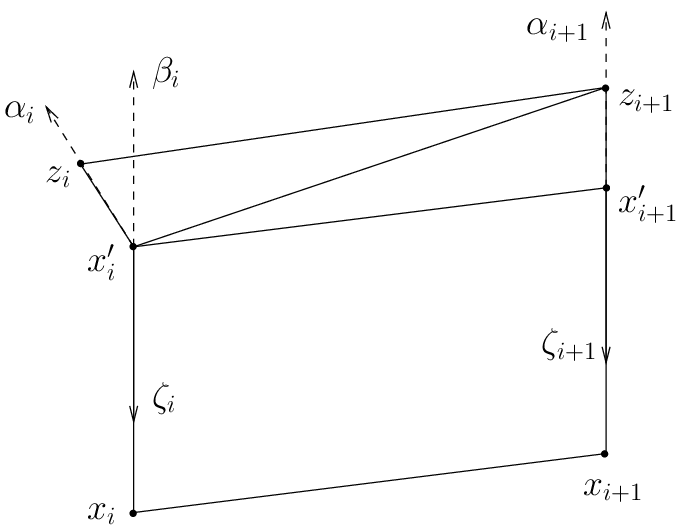}
\caption{Extending the geodesic segments $\nu_i$}\label{fig:localparallelset}
\end{center}
\end{figure}

If $\alpha_i=\beta_i$ for all $i$, then we can extend the $\nu_i$ inside
$A_{i-1}\cap A_{i}$ contradicting the maximality of $P$. Hence, there
is a $j$ such that $\alpha_j\neq \beta_j$.
Actually more is true:
if it holds for all $i$ that
$d(\ora{x_i'x_{i+1}'},\alpha_i)=d(\ora{x_i'x_{i+1}'},\beta_i)$,
then $\zeta_i\ora{x_i'x_{i+1}'}\alpha_i$ is a geodesic
segment in $\Sigma_{x_i'}X$ of length $\pi$.
Let $z_{i+1}\in A_{i}$ be a point near $x_{i+1}'$ with $\ora{x_{i+1}'z_{i+1}}=\alpha_{i+1}$
for $i=0,\dots,n-1$.
We can choose $z_{i+1}$ close enough to $x_{i+1}'$, so that
$\ora{x_i'z_{i+1}}$ is a regular point in the same Weyl chamber as $\ora{x_i'x_{i+1}'}$
because $p$ is a regular polygon.
It follows that $\ora{x_i'z_{i+1}}$ lies in the intersection of the segments
$\zeta_i\ora{x_i'x_{i+1}'}\alpha_i$ and $\zeta_i\ora{x_i'x_{i+1}'}\beta_i$.
Thus $\zeta_i\ora{x_i'z_{i+1}}\alpha_i$ is a geodesic segment of length $\pi$.
After perhaps taking $z_i\in A_i$ closer to $x_i'$
we can accomplish that
$CH(x_i',z_i,z_{i+1})$ is a flat triangle. 
It follows that the union of the (2-dimensional) flat convex sets 
$CH(x_i,x_{i+1},x_{i+1}',x_i')$, $CH(x_i',x_{i+1}',z_{i+1})$ and
$CH(x_i',z_{i+1},z_i)$ is a flat convex quadrilateral.
(See Figure~\ref{fig:localparallelset}.)
Notice also that the segments $\nu_i(s^-)z_i$ 
are extensions of the geodesic segments 
$\nu_i(s^-)\nu_i(s^+)$. Thus this contradicts as well the maximality of $P$.
Hence, there is a $j$ such that
$d(\ora{x_j'x_{j+1}'},\alpha_j) > d(\ora{x_j'x_{j+1}'},\beta_j)$.
Analogously, there is a $k$ such that 
$d(\ora{x_k'x_{k-1}'},\beta_k) > d(\ora{x_k'x_{k-1}'},\alpha_k)$.

For some small $\epsilon > 0$ let $\tilde x_j \in A_{j-1}$
be a point such that $d(x_j',\tilde x_j)=\epsilon$ and
$\ora{x_j'\tilde x_j} = \alpha_j$.
Then $\sigma(\tilde x_jx_{j+1}')
=\sigma(x_j'x_{j+1}')-\epsilon \cdot ov_{\tilde\eta}
=\sigma(x_jx_{j+1})-\epsilon\cdot ov_{\tilde\eta}$ for some vertex 
$\tilde\eta\in \partial_T E$ of the same type as $\eta_{j}$
(compare with Section~\ref{sec:functional}).
By the above consideration we must have
$\tilde\eta\neq\eta_{j}$, otherwise
$d(\ora{x_j'x_{j+1}'},\alpha_j) = d(\ora{x_j'x_{j+1}'},\beta_j)$
(recall that $l_{\eta_j}(\sigma(x_jx_{j+1}))=
\langle x_jx_{j+1},x_j\exp_{x_j}(1\!\cdot\!\xi_j)\rangle$ and 
$\dot\nu_j(0)=\ora{x_j\xi_j}$).
It follows that $l_{\eta_{j}}(\sigma(\tilde x_jx_{j+1}'))=
\langle \sigma(x_jx_{j+1})-\epsilon\cdot ov_{\tilde\eta}\,,\;ov_{\eta_{j}}\rangle
 = l_{\eta_{j}}(\sigma(x_jx_{j+1})) - \epsilon \langle ov_{\tilde\eta}, ov_{\eta_{j}}\rangle
 > l_{\eta_{j}}(\sigma(x_jx_{j+1})) - \epsilon$.
On the other hand,
$\sigma(x_{j-1}'\tilde x_j)=\sigma(x_{j-1}'x_j' + x_j'\tilde x_j) =
\sigma(x_{j-1}x_j) + \epsilon\cdot ov_{\eta_{j-1}}$	
and this implies
$l_{\eta_{j-1}}(\sigma(x_{j-1}'\tilde x_j))=
\langle \sigma(x_{j-1}x_j) + \epsilon\cdot ov_{\eta_{j-1}}\,,\;ov_{\eta_{j-1}}\rangle=
l_{\eta_j}(\sigma(x_{j-1}x_j))+\epsilon$.
Thus, for $p_1:=(x_1',\dots,\tilde x_j, \dots, x_{n-1}')$
we have $L(\sigma(p_1))>L(\sigma(p))=0$.

Let now $\tilde x_k \in A_{k}$ be a point such that 
$d(x_k',\tilde x_k)=\epsilon$ and
$\ora{x_k'\tilde x_k} = \beta_k$.
It follows
analogously for $p_2:=(x_1',\dots,\tilde x_k, \dots, x_{n-1}')$
that $L(\sigma(p_2))<L(\sigma(p))=0$.
\qed

The next question is what happens when $p$ is contained in such a parallel set
$P_{\xi,\xi'}$. In this last situation we cannot always get the same conclusion
as in Lemmata~\ref{lem:wallbyaccident} and \ref{lem:localparallelset}.
For instance, if the wall $H_L$ lies in the boundary of $\polspace$, then we can
cross $H_L$ in one direction but not in the opposite one.

\begin{rem}
Suppose that $p$ is contained in $P_{\xi,\xi'}$
with $\xi,\xi'$ as in Lemma~\ref{lem:localparallelset}. 
Let $b_{\xi'}:X\rightarrow\R$ be a Busemann function associated to 
$\xi'$ (see e.g. \cite[Sec. 2.2]{KLM2} for a definition). 
Then by considering an apartment in $P_{\xi,\xi'}$ containing the side
$x_{i}x_{i+1}$, we see that
$l_{\eta_i}(\sigma(x_{i}x_{i+1}))=b_{\xi'}(x_{i+1})-b_{\xi'}(x_{i})$.
In particular
$$
L(\sigma(p))=l_{\eta_0}(\sigma(x_0x_1))+\dots + l_{\eta_{n-1}}(\sigma(x_{n-1}x_0))
=\sum_{i=0}^{n-1} (b_{\xi'}(x_{i+1})-b_{\xi'}(x_{i}))=0.
$$

Thus, if $p'$ is the result of a variation
of the polygon $p$ within the parallel set
$P_{\xi,\xi'}$, it still holds $L(\sigma(p'))=0$.
\end{rem}

The next lemma gives a condition that let us cross the wall $H_L$ in the
positive direction. 

\begin{figure}[h]
\begin{center}
 \includegraphics[scale=0.6]{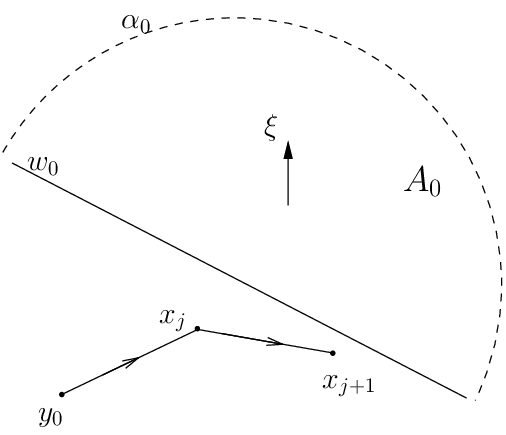}
\caption{Setting of Lemma~\ref{lem:crossablewalls}}\label{fig:samehalfspace}
\end{center}
\end{figure}

Suppose $p$ is contained in $P_{\xi,\xi'}$ where 
$\xi,\xi'\in \partial_T X$ are antipodal points
such that 
$\ora{x_i\xi_i}=\ora{x_i\xi}$ for all $i$. 
Assume also that there are vertices $x_i,x_j,x_{j+1}$
of $p$ with the following property.
Let $A_0, A_1$ be apartments in $P_{\xi,\xi'}$ containing the segment
$x_jx_{j+1}$ and an initial part of the segment $x_jx_i$ and 
$x_{j+1}x_i$ respectively. 
Let $y_k\in A_k$ for $k=0,1$ be points in the initial parts
of the segments
$x_jx_i$ and $x_{j+1}x_i$ respectively. Thus $x_jx_{j+1}y_k$ are flat triangles
in $A_k$.
Suppose that for some $k=0,1$ there is a root $\alpha_k\subset \partial_T A_k$
such that the directions $\xi$, $\ora{x_jx_{j+1}}$
and $(-1)^k\ora{y_kx_{j+k}}$ 
lie in the interior of $\alpha_k$ (after the 
natural identification of $\partial_T A_k$ and $\Sigma_x A_k$ for $x\in A_k$).
By $-\ora{y_1x_{j+1}}$ we mean $\ora{x_{j+1}y_1}$.
(See Figure~\ref{fig:samehalfspace}.)

\begin{lem}\label{lem:crossablewalls}
Under the assumptions above, for any neighborhood
$U$ of $\sigma(p)$ in $\eucch^n$ there is an $n$-gon $p'$ in $X$
with $\sigma(p')\in U$ and $L\circ\sigma(p')>0$.
\end{lem}
\proof
We show the lemma when the root $\alpha_k\subset \partial_T A_k$ exists for $k=0$. 
The other case $k=1$ is analogous.

Let $w_0\subset A_0$ be a wall, such that $\alpha_0$ is bounded by $\partial_T w_0$.
Let $h^\pm\subset A_0$ be the roots bounded by $w_0$ with ideal boundary $\alpha_0$ and
its antipodal, respectively.
Let $\epsilon >0$ be small.
First we displace the polygon $p$ in $P_{\xi,\xi'}$ 
along $\xi\xi'$
such that $x_j$ lies in
$h^+$ and $d(x_j,w_0)<\epsilon$.
By the hypothesis we can take $\epsilon$ small enough so that $y_0\in h^-$.
Let $A_0'$ be an apartment in $X$ such that $A_0\cap A_0' = \overline{h^-}$.
Let $x_j'\in A_0'$ be the point such that $d(y_0,x_j')=d(y_0,x_j)$
and $\ora{y_0x_j'}=\ora{y_0x_j}$.
Let $z\in A_0$ be the reflection of $x_j$ in the hyperplane $w_0$.

Observe that $x_j'\notin A_0$ and $x_{j+1}\notin A_0'$.
It follows that $\sigma(x_j'x_{j+1})=\sigma(zx_{j+1})$.
In particular 
\begin{eqnarray*}
 l_{\eta_{j}}(\sigma(x_j'x_{j+1})) & = & l_{\eta_{j}}(\sigma(zx_{j+1}))
    = \langle zx_j+x_jx_{j+1},\, z\exp_z(1\!\cdot\!\xi) \rangle \\
& = & \langle zx_j,\,z\exp_z(1\!\cdot\!\xi)\rangle + l_{\eta_{j}}(\sigma(x_jx_{j+1}))
    \,>\,l_{\eta_{j}}(\sigma(x_jx_{j+1})).
\end{eqnarray*}

Notice that the refined length of $x_j'x_i$ is the same as of $x_jx_i$. 
Hence, by Theorem~\ref{thm:transfer}
 we can transfer the polygon $(x_i,x_{i+1},\dots,x_j)$ to a polygon
$(x_i',x_{i+1}',\dots,x_j')$ with $x_i'=x_i$ and with 
the same $\Delta$-valued side lengths.
The $n$-gon $p' = (x_i',x_{i+1}',\dots,x_j',x_{j+1},\dots,x_{i-1})$
satisfies the conclusion of the lemma.
\qed

If the polygon $p$ is completely contained in an apartment in $P_{\xi,\xi'}$,
then the condition for the lemma above can be stated more easily.

\begin{cor}\label{cor:polygoninapt}
Suppose $p$ is contained in an apartment $A\subset P_{\xi,\xi'}$.
Suppose that there are two sides $x_ix_{i+1}$, $x_jx_{j+1}$ of $p$
and a root $\alpha\subset A$, such that
the directions $\xi$, $\ora{x_ix_{i+1}}$
and $\ora{x_jx_{j+1}}$ 
lie in the interior of $\alpha$.
Then for any neighborhood
$U$ of $\sigma(p)$ in $\eucch^n$ there is an $n$-gon $p'$ in $X$
with $\sigma(p')\in U$ and $L\circ\sigma(p')>0$.
\end{cor}
\proof
Consider the oriented segments $d_1=x_ix_j$ and $d_2=x_jx_i$. 
After a small variation
of the polygon $p$ inside of the apartment $A$, we may assume that $d_1$
(and therefore also $d_2$) is regular.
Then for one $k=1,2$, it must hold, that $\ora{d_k}$ 
lies in the interior of $\alpha$. If $k=1$,
then Lemma~\ref{lem:crossablewalls} applies for 
the vertices $x_i,x_j,x_{j+1}$ and if $k=2$, then it applies
for the vertices $x_j,x_i,x_{i+1}$.
\qed

Let us assume now that the building $X$ has rank 2. 
We explain another method special for this case
to cross the wall $H_L$.

\begin{figure}[h]
\begin{center}
 \includegraphics[scale=0.6]{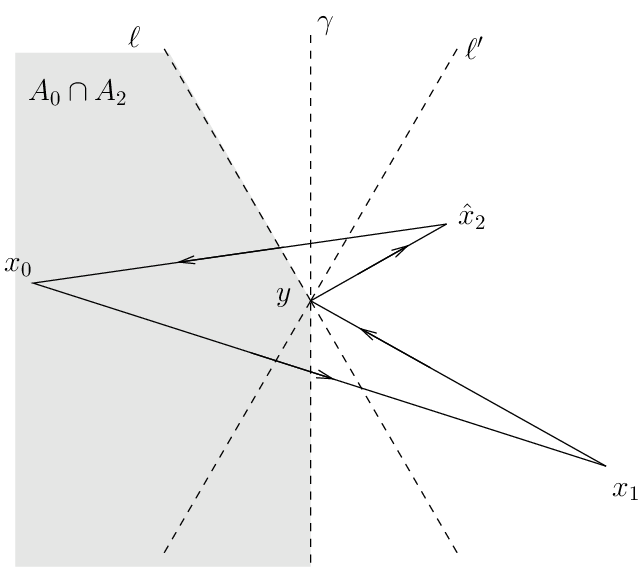}
\caption{The folded triangle $\hat p$}\label{fig:variationrank2a}
\end{center}
\end{figure}

Let $p=(x_0,x_1,x_2)$ be a regular triangle contained in $P_{\xi,\xi'}$ 
but not contained in any apartment.
It is easy to see, that when we fold $p$ into an apartment $A$, 
it has exactly one break point.
After relabeling the vertices we can assume that the break point $y$ 
lies in the 
side $x_1x_2$ and that the sides of the folded triangle 
$\hat p = (\hat x_0=x_0, \hat x_1=x_1, y, \hat x_2)$ 
do not intersect in their interiors (see Figure~\ref{fig:variationrank2a}).
After displacing $\hat p$ in $P_{\xi,\xi'}$
along $\xi\xi'$ we can assume that $y$ is a vertex
of $X$. 
Let $\gamma\subset A$ be the singular line through $y$ connecting $\xi$ and $\xi'$.

\begin{lem}\label{lem:crossablewallsrank2}
We use the setting above (in particular, $rank(X)=2$).
Suppose that the Weyl chamber containing $\ora{yx_1}$
is not adjacent to $\Sigma_y\gamma$.
Then for any neighborhood
$U$ of $\sigma(p)$ in $\eucch^3$ there are triangles
$p_1,p_2$ in $X$
with $\sigma(p_1),\sigma(p_2)\in U$ and $L\circ\sigma(p_1)>0>L\circ\sigma(p_2)$.
\end{lem}
\proof
Let $\ell\subset A$ be the singular line through $y$ such that 
$\Sigma_y\ell$ is adjacent to
the simplicial convex hull of $\ora{yx_1}\ora{y\hat x_2}$
and the directions $\xi$, $\ora{yx_1}$
and $\ora{y\hat x_2}$ are contained in the interior of the same 
root bounded by $\partial_T\ell$. 
It exists by the assumptions of the lemma.
Let $\ell'\subset A$ be the reflection of $\ell$ in $\gamma$.
Let $h_\gamma^-\subset A$ be the root bounded by $\gamma$ containing $x_0$
and let $h_\gamma^+\subset A$ be the antipodal root.
Similarly, let $h_\ell^+,h_{\ell'}^+\subset A$ be the roots bounded by $\ell,\ell'$
containing $\xi$ and let $h_\ell^-,h_{\ell'}^-\subset A$ be the antipodal roots.
Then the simplicial convex hull
of $\ora{yx_1}\ora{y\hat x_2}$ is 
$\Sigma_y(\overline{h_\ell^+}\cap \overline{h_{\ell'}^-})$. 
(See Figure~\ref{fig:variationrank2a}.)

Let $A_2$ be an apartment in $X$ such that 
$A\cap A_2=\overline{h_\gamma^-}\cap\overline{h_\ell^-}$.
Let $x_2'\in A_2$ be the point so that 
$d(x_0,x_2')=d(x_0,\hat x_2)$ and $\ora{x_0x_2'}=\ora{x_0\hat x_2}$.
Notice that $x_2'\notin A$, thus, $x_2'\neq\hat x_2$. Observe also that
$x_2'\notin P_{\xi,\xi'}$, hence, $x_2'\neq x_2$.

\begin{figure}[h]
\begin{center}
 \includegraphics[scale=0.6]{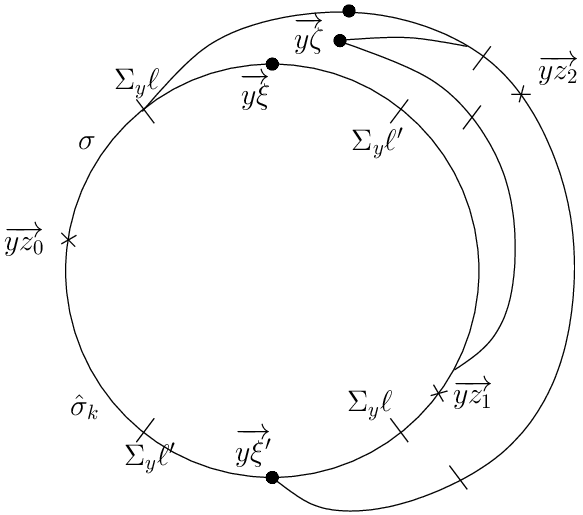}
\caption{$\Sigma_y X$}\label{fig:variationrank2b}
\end{center}
\end{figure}

The concatenation of the segments $\ora{yx_1}\ora{y\xi'}\in\Sigma_y A$ and 
$\ora{y\xi'}\ora{yx_2'}\in\Sigma_y A_2$ gives a segment in $\Sigma_y X$
of length $\pi$ (see Figure~\ref{fig:variationrank2b}). 
Therefore $x_1yx_2$ is a geodesic segment
and the triangle $p'=(x_0,x_1,x_2')=:(z_0,z_1,z_2)$ has the same side lengths as $p$.
Set $A_0:=A$ and
let $A_1$ be an apartment in $X$ containing the segment $z_1z_2$.

Let $\nu_i$ be the geodesic rays with $\nu_i(0)=z_i$ and $\nu_i(-\infty)=\xi'$.
Then $CH(\nu_i,\nu_{i+1})$ are (2-dimensional) flat stripes.
We want to see that the $\nu_i$ cannot be extended to parallel geodesic lines.
Suppose then the contrary: there are parallel geodesic lines $\nu_i'$ containing
$\nu_i$. Set $\zeta:=\nu_i'(\infty)$.
Then $p'\subset Y:=P_{\zeta,\xi'}$ and in particular,
$\ora{y\xi'},\ora{yz_i}\in \Sigma_y Y$.
Since $\ora{yz_1},\ora{yz_2}\in\Sigma_yA_1$ are antipodal regular points, the 
apartment containing them is unique. 
Therefore $\Sigma_yA_1\subset \Sigma_y Y$ and in particular,
$\ora{y\zeta}\in\Sigma_yA_1$.

Let $k\in\{1,2\}$ be so that the Weyl chamber 
containing $\ora{y\hat x_k}$ is adjacent to $\Sigma_y\ell'$.
Let $\sigma_k\subset \Sigma_yA_{2k-2}$
be the Weyl chamber containing $\ora{yz_k}$ and let 
$\hat\sigma_k\in\Sigma_y(A_0\cap A_2)$ be the antipodal chamber to $\sigma_k$.
(See Figure~\ref{fig:variationrank2b} for $k=2$.)
Notice that $\ora{y\xi'}\ora{yz_0}$ intersects $\hat\sigma_k$ in its interior.
In particular $\hat\sigma_k\subset \Sigma_y Y$. 
It follows that the unique apartment containing 
$\sigma_k$ and $\hat\sigma_k$ is contained in $\Sigma_y Y$, i.e.\
$\Sigma_y A_{2k-2}\subset \Sigma_y Y$.

Let $\sigma\subset\Sigma_y(A_0\cap A_2)
\subset \Sigma_y A_{2k-2}\subset\Sigma_y Y$
 be the Weyl chamber adjacent to
$\ell$. The Weyl chamber containing $\ora{yz_{3-k}}$ is antipodal to $\sigma$.
Hence, the unique apartment containing $\sigma$ and $\ora{yz_{3-k}}$
is contained in $\Sigma_y Y$, i.e.\ $\Sigma_y A_{4-2k}\subset \Sigma_y Y$.

We have conclude that $A_0,A_2\subset \Sigma_y Y = \Sigma_y P_{\zeta,\xi'}$,
but this is not possible because of the construction of $A_2$.
Therefore the geodesic rays $\nu_i$ cannot be extended to complete parallel
 geodesic lines. The lemma now follows from Lemma~\ref{lem:localparallelset}
and its proof.
\qed

We can show now that for rank 2 the space $\polspace$ is a polyhedral cone. 
Its convexity will be shown in the next section.

\begin{prop}\label{prop:unionofcones}
If $X$ has rank 2, then $\polspace$ is a union of the closures of 
polyhedral cones in ${\cal C}_n$.
\end{prop}
\proof
We have already seen in Corollary~\ref{cor:opencones} that if 
for $C\in{\cal C}_n$ holds $\polspace\cap C\neq\emptyset$, then 
$\overline C \subset \polspace$. 
Now let $p=(x_0,\dots,x_{n-1})$ be a polygon in $X$. 
We want to show that $\sigma(p)$ is
contained in $\overline C$ for some $C\in{\cal C}_n$
with $C\subset \polspace$. 
Since any polygon
can be approximated by regular polygons, we may assume that $p$ is regular.
Suppose now $s:=\sigma(p)\in H_L$ with $L=(l_{\eta_0},\dots,l_{\eta_{n-1}})$. 
If for any neighborhood $U$ of $s$ we can find polygons with side lengths in 
$U\setminus H_L$, then we are done. Indeed, in this case, there is an open cone
$C\in{\cal C}_n$ such that $\polspace\cap C\neq\emptyset$ and $s\in \overline C$.

Suppose then that for some neighborhood $U$ of $\sigma(p)$
we cannot find polygons $p'$ with side lengths in $U$
and $L\circ\sigma(p')\neq 0$.
Lemmata~\ref{lem:wallbyaccident} and \ref{lem:localparallelset}
implies that $p$ lies in a parallel set $P_{\xi,\xi'}$ and the functional
$L$ is given in $p$ by taking scalar product with a unit vector in 
the direction of $\xi$.
Suppose first that the triangle $t=(x_0,x_1,x_2)$ 
lies completely in an apartment in $P_{\xi,\xi'}$. 
Then it is easy to see that Corollary~\ref{cor:polygoninapt}
must apply for one of the functionals 
$L':=(l_{\eta_0},l_{\eta_1},l_{\eta'})$ or $-L'$,
where $\eta'$ is so that $l_{\eta'}(\sigma(x_2x_0))=
\langle x_2x_0, \,x_2\!\exp_{x_2}(1\!\cdot\!\xi)\rangle$.
If $t$ is not contained in an apartment, then
we fold it into an apartment as in the setting of Lemma~\ref{lem:crossablewallsrank2}. 
Then, either Lemma~\ref{lem:crossablewallsrank2} applies or
the Weyl chamber containing
the direction $\ora{x_ix_{i+1}}$ of the side of $t$ with the break point
must be adjacent to $\ora{x_i\xi}$ or $\ora{x_i\xi'}$. 
If the last occurs, it is again easy to see, that
Lemma~\ref{lem:crossablewalls}
must apply for $L'$ or $-L'$.
In any case, 
we find a triangle $t'=(x_0',x_1',x_2')$ with $L'(\sigma(t))\neq 0$
and such that (modulo displacement in $P_{\xi,\xi'}$ along $\xi\xi'$) 
the refined side lengths of
$t'$ are as near as we want to the ones of $t$.
After a small variation of the polygon $(x_0,x_2,\dots,x_{n-1})$ inside the 
parallel set $P_{\xi,\xi'}$ and displacing it along $\xi\xi'$, we obtain
a polygon $q=(x_0'',x_2'',\dots,x_{n-1}'')$ so that the refined side length of
$x_0'x_2'$ is the same as of $x_0''x_2''$.
Since $q$ is still contained in the parallel set $P_{\xi,\xi'}$, it holds
$L''(\sigma(q))=0$ for $L''=(l_{\eta''},l_{\eta_2},\dots,l_{\eta_{n-1}})$
where $\eta''$ is so that $l_{\eta''}(\sigma(x_0''x_2''))=
\langle x_0''x_2'', \,x_0''\!\exp_{x_0''}(1\!\cdot\!\xi)\rangle$.
In particular, $l_{\eta'}(\sigma(x_2'x_0'))=-l_{\eta''}(\sigma(x_0''x_2''))$.
Then by the Transfer Theorem~\ref{thm:transfer} we can glue $t'$ and $q$
along $x_0'x_2'$ and $x_0''x_2''$ to a polygon $p'$ with $\Delta$-valued
side lengths near $s$ and 
$L(\sigma(p')) = L'(\sigma(t'))+L''(\sigma(q))\neq 0$.
\qed

\begin{rem}
 Proposition~\ref{prop:unionofcones} is also true in rank $>2$ by the
results of \cite{KLM1} and \cite{KLM2}. However our proof here uses 
Lemma~\ref{lem:crossablewallsrank2}, which we only showed in rank 2.
\end{rem}

\subsection{The boundary of $\polspace$}\label{sec:boundary}

We have seen in the previous section different methods which allows to 
cross certain walls $H_L$ within the space $\polspace$.
We will show in this section that 
for the case of buildings of rank 2
the walls where this method cannot be
applied are precisely the walls that determine the boundary 
of $\polspace$. 
That is, if a wall cannot be crossed with the polygon variations
of Section~\ref{sec:crosswalls}, it is because that wall cannot be crossed
at all.

First we characterize the walls $H_L$ that cannot be crossed 
with the methods above
in terms of the
combinatorics of the associated spherical Coxeter complex $(S,W)$.
We use the same notation as in Section~\ref{sec:functional}.
Let $\Delta:=\partial_T\eucch \subset \partial_T E=S$ 
be the spherical Weyl chamber. 
Let $\delta\in W$ be the element of the Weyl group such that $\Delta$
and $\delta\Delta$ are antipodal. Notice that $\delta^2=Id$.
We say that a root $\alpha\subset S$ is a {\em positive root} if 
$\Delta\subset \alpha$. Let $\Lambda^+$ denote the set of positive roots of $(S,W)$.
A root which is not a positive root is called a {\em negative root}.

For each element $\omega\in W$ we define a subset of $\Lambda^+$:
$$
T_\omega:=\{\alpha\in \Lambda^+ \;|\; \omega^{-1}\Delta\subset\alpha\}
=\Lambda^+ \cap \omega^{-1}\Lambda^+.
$$

Let $\eta\in\Delta$ be a vertex of the spherical Weyl chamber
and let $\Lambda_\eta^+$ be the set of positive roots containing
$\eta$ in their boundaries.
Let $B_n(\eta)$ be the set of $n$-tuples 
$(\omega_0,\dots,\omega_{n-1})\in W^n$ such that
\begin{enumerate}
\item[($\ast$)] $T_{\omega_i}\cap T_{\omega_j}
	\subset \Lambda_\eta^+$ for all $i\neq j$,
\item[($\ast\ast$)] $\bigcup\limits_{i=0}^{n-1} T_{\omega_i}= \Lambda^+$.
\end{enumerate}

For $\bar\omega=(\omega_0,\dots,\omega_{n-1})\in W^n$ we write
$\bar\omega\eta:=(\omega_0\eta,\dots,\omega_{n-1}\eta)\in(W\eta)^n$.
and for $\bar\eta=(\eta_0,\dots,\eta_{n-1})\in (W\eta)^n$ set
$L_{\bar\eta}=(l_{\eta_0},\dots,l_{\eta_{n-1}})$.
Let ${\cal B}_n\subset {\cal L}_n$ be the union of the sets 
$\{L_{\bar\eta} \;|\; \bar\eta=\bar\omega\eta,\: \bar\omega\in B_n(\eta)\}$ for all
vertices $\eta\in \Delta$. 

We will see in Lemma~\ref{lem:notpassablewalls} below that
the walls $H_L$ that cannot be crossed in the positive direction
with our previous methods
are precisely the ones of the form $L_{\bar\eta}$
with $\bar\eta=\bar\omega\eta$ 
and $\bar\omega=(\omega_0,\dots,\omega_{n-1})$ satisfying the property ($\ast$).
A motivation for this property ($\ast$) can already be seen in 
Corollary~\ref{cor:polygoninapt}.
Namely, let $p=(x_0,\dots, x_{n-1})$ 
be a regular polygon contained in an apartment $A$ and 
let $\eta'$ be a vertex of $\partial_T A$. Let $L$ be the functional that
in $p$ corresponds to taking scalar product with a unit vector
in the direction of $\eta'$. In particular, $L(\sigma(p))=0$.
Fix a Weyl chamber $\tau\subset \partial_T A$
containing $\eta'$ and let $\tau_i\subset \partial_T A$ be the Weyl chamber
containing the direction $\ora{x_ix_{i+1}}$ after identifying 
$\Sigma_{x_i}A\cong\partial_T A$.
Let $\omega_i\in W$ be such that $\tau=\omega_i\tau_i$.
Let $\eta\in\Delta$ be the vertex of the same type as $\eta'$.
Then $L=(l_{\omega_0\eta},\dots,l_{\omega_{n-1}\eta})$.
If
$T_{\omega_i}\cap T_{\omega_j}\not\subset \Lambda_\eta^+$ for some $i\neq j$,
then Corollary~\ref{cor:polygoninapt} applies and we can 
increase the functional $L$ near $\sigma(p)$.
The property ($\ast\ast$) is introduced to avoid later 
obvious redundancies in the set of generalized triangle inequalities.
This can be seen in the Lemma~\ref{lem:noirredundant}.

\begin{lem}\label{lem:ineq=weakineq}
Let $\bar\eta=(\eta_0,\dots,\eta_{n-1})\in (W\eta)^n$.
If there exist $j,j'\in\{0,\dots,n-1\}$, $j\neq j'$ with
$\eta_{j'}=\delta\eta_{j}$
and $\delta\eta_i=\eta\in\Delta$ for $i\neq j,j'$, then
$\bar\eta\in B_n(\eta)\!\cdot\!\eta$.
If $(E,W)$ has rank 2, then the converse is also true.
\end{lem}
\proof To prove the first statement
choose $\omega_j\in W$ with $\omega_j\eta=\eta_j$
and set $\omega_{j'}:=\delta\omega_j$. 
Then $\omega_{j'}\eta=\delta\omega_j\eta=\delta\eta_j=\eta_{j'}$.
It follows that
$\omega_{j'}^{-1}\Delta$ is antipodal to $\omega_{j}^{-1}\Delta$
and $T_{\omega_j} \cap  T_{\omega_{j'}} = \emptyset$,
$T_{\omega_j} \cup  T_{\omega_{j'}}=\Lambda^+$.
Let $\omega_i=\delta$ for $i\neq j,j'$, then
$\eta_i=\omega_i\eta$ and $T_{\omega_i}=\emptyset$.
Hence, $\bar\omega=(\omega_0,\dots,\omega_{n-1})\in B_n(\eta)$ and
$\bar\eta = \bar\omega\eta$.	

\begin{figure}[h]
\begin{center}
 \includegraphics[scale=0.5]{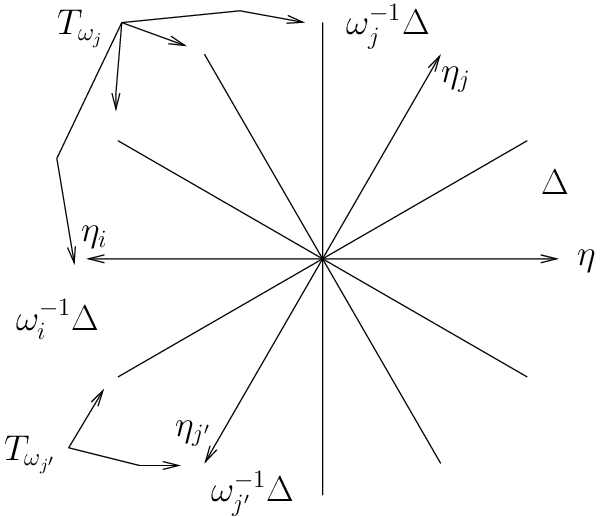}
\caption{${\cal B}_n^w$: weak triangle inequalities}\label{fig:Bnrank2}
\end{center}
\end{figure}

Now suppose that $(E,W)$ has rank 2 and let
$\bar\omega\in B_n(\eta)$ with $\bar\eta=\bar\omega\eta$. 
By property ($\ast\ast$) there is a $j$
with $T_{\omega_j}\not\subset\Lambda_\eta^+$,
that is, $-\eta\notin \omega_j^{-1}\Delta$. 
If $\eta\in\omega_j^{-1}\Delta$ (or equivalently, $\omega_j\in Stab_W(\eta)$)
then $\eta_j=\omega_j\eta=\eta$ and $\Lambda^+ \setminus T_{\omega_j}\subset \Lambda_\eta^+$.
This and properties ($\ast$) and ($\ast\ast$) imply
$T_{\omega_i}\subset \Lambda_\eta^+$ for $i\neq j$.
It follows that 
$-\eta\in\omega_i^{-1}\Delta$, or equivalently, $\delta\omega_i\in Stab_W(\eta)$.
Hence, $\eta_i=\omega_i\eta=\delta\delta\omega_i\eta=\delta\eta$
for $i\neq j$ and the assertion follows.
Suppose now that $-\eta,\eta\notin\omega_j^{-1}\Delta$.
Then there are two positive roots $\alpha_1,\alpha_2\notin\Lambda_\eta^+$
with $\omega_j^{-1}\Delta=\alpha_1\cap-\alpha_2$.
Then $\alpha_2\in \Lambda^+ \setminus T_{\omega_j}$. 
Let $j'$ be such that $\alpha_2\in T_{\omega_{j'}}$.
By property ($\ast$) it follows that $j'$ is unique and
$\alpha_1\notin T_{\omega_{j'}}$.
Hence, $\omega_{j'}^{-1}\Delta=-\alpha_1\cap\alpha_2$ 
and $\omega_{j'}^{-1}\Delta$
is antipodal to $\omega_j^{-1}\Delta$.
This implies that
$\omega_{j'}=\delta\omega_j$ and $\eta_{j}=\delta\eta_{j'}$.
Since $\Lambda^+ = T_{\omega_j} \cup  T_{\omega_{j'}}$,
by property $(\ast)$ we deduce that $T_{\omega_i}\subset\Lambda_\eta^+$ for $i\neq j,j'$
and this in turn implies that $\eta_i=\delta\eta$ for $i\neq j,j'$ (Fig.~\ref{fig:Bnrank2}).
\qed

\begin{rem}\label{rem:weakineq}
Let ${\cal B}_n^w\subset {\cal L}_n$ be the set 
of functionals $L_{\bar\eta}$ for $\bar\eta$ satisfying
the hypothesis of the previous Lemma.
The inequalities $L\leq 0$ for $L\in {\cal B}_n^w$ 
are the so-called 
{\em weak triangle inequalities} (cf. \cite[Section 3.8]{KLM1}).
Thus, Lemma~\ref{lem:ineq=weakineq} just states that 
${\cal B}_n^w\subset {\cal B}_n$ and that for rank 2 it actually holds
${\cal B}_n={\cal B}_n^w$ (Fig.~\ref{fig:Bnrank2}).
\end{rem}

\begin{thm}[Weak Triangle Inequalities, {\cite[Theorem 3.34]{KLM1}}]
\label{thm:weakineq}
For any $n$-gon $p$ in $X$ and any functional $L\in {\cal B}_n^w$ holds
$L\circ\sigma(p)\leq 0$. That is,
$$
\polspace \subset \bigcap_{L\in {\cal B}_n^w} \{L\leq 0\}.
$$
\end{thm}
\proof
Let $p=(x_0,\dots,x_{n-1})$ be an $n$-gon in $X$. 
For the functional $L=(l_{\eta_0},\dots,l_{\eta_{n-1}})\in {\cal B}_n^w$,
let $j,j'$ be such that $\eta_{j}=\delta\eta_{j'}$
and $\delta\eta_i=\eta\in\Delta$ for $i\neq j,j'$.
Notice that 
$l_{\delta\eta}\leq l_{\eta'}$ in $\eucch$ for all $\eta'$ of the same type as $\eta$.
That is, $l_{\eta_i}$ is the smallest functional of the same type as $\eta$.
After shifting the subindices of the polygon and the functional
we can assume that $j'=0$. Let $\omega_0\in W$ be so that
$\omega_0\eta=\eta_0$.

Suppose first that $j=n-1$.
Fold the polygon $p$ into an apartment $A$, so that the broken sides are
$x_1x_2,\dots,x_{n-2}x_{n-1}$.
Let $\rho:A\rightarrow E$ be an isometry that sends $x_0$ to the vertex
$o$ of $\eucch\subset E$, induces an isomorphism of the Coxeter complexes
$(\partial_T A, W)$ and $(E,W)$ and so that $\rho(x_0x_1)\subset\omega_0^{-1}\eucch$.
Notice that $\rho$ is not necessarily
an isomorphism of Euclidean Coxeter complexes with the affine Weyl group
$W_{aff}$. Denote with $q$ the image under $\rho$ of the folded polygon.
By folding $E$ onto the Euclidean Weyl chamber $\omega_0^{-1}\eucch$ with the natural
``accordion'' map,
 we obtain a further folded polygon 
$q'=(y_0,\dots,y_k)$ where $y_0$ is the vertex of $\eucch$
and the $\Delta$-valued side lengths of the segments
$y_0y_1, y_0y_k \subset \omega_0^{-1}\eucch$ 
are the same as for $x_0x_1$ and $x_0x_{n-1}$ respectively.
Observe that $q'$ is not necessarily a billiard polygon in $(E,W_{aff})$, 
but if the side
$x_rx_{r+1}$ of $p$ is broken in $q'$ to the sides $y_sy_{s+1},y_{s+1}y_{s+2},
\dots,y_{t-1}y_t$, then the vectors
$\sigma(y_sy_{s+1}),\dots,\sigma(y_{t-1}y_t)$ are just multiples of 
$\sigma(x_rx_{r+1})$. 
This means, that if $W_{aff}'$ is the group generated by $W_{aff}$
and the whole translation group of $E$, then $q'$ is a billiard
polygon in $(E,W_{aff}')$. 
Notice also that for $i\neq 0,n-1$ holds
$l_{\eta_i}(\sigma(y_ly_{l+1}))\leq 
\langle y_ly_{l+1}, ov_{\eta} \rangle$
because of the observation at the beginning of the proof.
It follows that
$$
l_{\eta_1}(\sigma(x_1x_2))+\dots+l_{\eta_{n-2}}(\sigma(x_{n-2}x_{n-1})) \leq
\langle y_1y_2, ov_{\eta} \rangle + \dots + \langle y_{k-1}y_{k}, ov_{\eta} \rangle =
\langle y_1y_{k}, ov_{\eta} \rangle.
$$

On the other hand, since $y_0y_1, y_0y_k \subset \omega_0^{-1}\eucch$ and
$\eta_{n-1}=\delta\eta_0=\delta\omega_0\eta$ it follows that
$l_{\eta_0}(\sigma(x_0x_1))=l_{\eta_0}(\sigma(y_0y_1))
= \langle \sigma(y_0y_1), ov_{\eta_0} \rangle=\langle y_0y_1, ov_{\eta} \rangle$ 
and
$l_{\eta_{n-1}}(\sigma(x_{n-1}x_0))=l_{\eta_{n-1}}(\sigma(y_ky_0))=
\langle \sigma(y_ky_0), ov_{\eta_{n-1}} \rangle = \langle y_ky_0, ov_{\eta} \rangle$.
Hence,
$L(\sigma(p))\leq \langle y_0y_1, ov_\eta \rangle
+\langle y_1y_k, ov_\eta \rangle
+\langle y_ky_0, ov_\eta \rangle = 0$.

The general case (i.e.\ $j\in\{1,\dots,n-1\}$)
now follows from the special case above by 
considering the polygons 
$p_1=(x_0,x_1,\dots,x_j)$, 
and
$p_2=(x_0,x_j,x_{j+1},\dots,x_{n-1})$ 
with the functionals 
$L_1=(l_{\eta_0},l_{\eta_1},\dots,l_{\eta_{j}})$
respectively
$L_2=(l_{\eta_{0}},l_{\eta_j},l_{\eta_{j+1}},\dots,l_{\eta_{n-1}})$.
Indeed, notice that since $\eta_0=\delta\eta_j$, it holds
$l_{\eta_{j}}(\sigma(x_{j}x_0))= \langle \sigma(x_jx_0), ov_{\eta_{j}} \rangle =
\langle -\sigma(x_0x_j), ov_{\eta_{0}} \rangle =
-l_{\eta_{0}}(\sigma(x_0x_{j}))$.
Hence, $L(\sigma(p)) = L_1(\sigma(p_1)) + L_2(\sigma(p_2))\leq 0$.
\qed

\begin{lem}\label{lem:noirredundant}
If $X$ has rank 2 and
$\bar\omega\in(W)^n$ satisfies the property ($\ast$) 
but not the property
($\ast\ast$) for some vertex $\eta\in\Delta$, 
then there is a $\bar\omega'\in B_n(\eta)$ so that
$L_{\bar\omega\eta}\circ\sigma(p)\leq L_{\bar\omega'\eta}\circ\sigma(p)$ 
for all $n$-gons $p$ in $X$. 
If $p$ is regular, then the strict inequality holds,
in particular, $L_{\bar\omega\eta}\circ\sigma(p) < 0$. 
\end{lem}
\proof
Let $\bar\omega=(\omega_0,\dots,\omega_{n-1})$ satisfy the property
($\ast$). It is easy to see that in rank 2 at most for two indices $i$
can hold $T_{\omega_i}\not\subset\Lambda_\eta^+$. 
Then we can find $j\neq j'$ so that
$T_{\omega_i}\subset\Lambda_\eta^+$ for all $i\neq j,j'$.
Let $\hat\omega_j:=\delta\omega_j$, then $\hat\omega_j^{-1}\eucch$ is antipodal
to $\omega_{j}^{-1}\eucch$.
If $\bar\omega$ does not satisfy the property ($\ast\ast$), then
$\hat\omega_j\neq\omega_{j'}$, moreover,
$T_{\omega_{j'}}\setminus \Lambda_\eta^+\subsetneq T_{\hat\omega_{j}} =
\Lambda^+ \setminus T_{\omega_{j}}$.
This implies that 
$l_{\omega_{j'}\eta}\leq l_{\hat\omega_j\eta}$ in 
$\eucch$ and since
$\hat\omega_j\neq\omega_{j'}$ the strict inequality holds for regular segments.
Thus we obtain $\bar\omega'\in B_n(\eta)$ by replacing $\omega_{j'}$ by
$\hat\omega_j$ in $\bar\omega$.
\qed

\begin{lem}\label{lem:notpassablewalls}
Suppose $X$ has rank 2 and let $p$ be a regular $n$-gon in $X$.
Suppose that $\sigma(p)\in H_L$ for some functional $L$
with $L,-L\in {\cal L}_n\setminus {\cal B}_n$.
Then for any neighborhood
$U$ of $\sigma(p)$ in $\eucch^n$ 
there exist $n$-gons $p_1,p_2$ in $X$
with $\sigma(p_i)\in U$ and $L\circ\sigma(p_1)>0>L\circ\sigma(p_2)$.
\end{lem}
\proof
Suppose that for a neighborhood $U$ of $\sigma(p)$ in $\eucch^n$,
we cannot find a polygon $p_1$ in $X$ with $\sigma(p_1)\in U$ and $L\circ\sigma(p_1)>0$. 
(The other inequality follows considering the functional $-L$.)
It follows from
Lemmata~\ref{lem:wallbyaccident} and \ref{lem:localparallelset} that 
$p$ lies in a parallel set $P_{\xi,\xi'}$ and the functional $L$
in $p$
is just given by taking scalar product with a unit vector in 
the direction of $\xi$.
Fold the polygon in an apartment $A\subset P_{\xi,\xi'}$ 
so that the broken sides are
$x_1x_2,\dots,x_{n-2}x_{n-1}$.
Let $\rho:A\rightarrow E$ be an isomorphism 
of Euclidean Coxeter complexes
that sends $\xi$ to the vertex $\eta\in\eucch$ of the same type. 

Suppose $X$ has only one vertex $o$ and let $\gamma\subset P_{\xi,\xi'}$ 
be the line through $o$ and $\gamma(\infty)=\xi,\gamma(-\infty)=\xi'$.
Then the break points of the folded polygon all lie on $\gamma$.
We may assume that the folded polygon has at most one break point because
any two consecutive break points can be simultaneously {\em unfolded}.
Let $k$ be so that the break point $y$ lies on the side $x_kx_{k+1}$
(if there is no break point we take $k=n-1$).
Then the folded polygon has the form 
$p'=(x_0,x_1,\dots,x_k,y,\hat x_{k+1},\dots,\hat x_{n-1})$.
Let $\omega_i\in W$ be so that 
$\omega_i^{-1}\Delta$ contains the direction $\rho(\ora{x_{i}x_{i+1}})$
for $0\leq i\leq k-1$, 
$\rho(\ora{x_ky})$ for $i=k$,
$\rho(\ora{\hat x_{i}\hat x_{i+1}})$
for $k+1\leq i \leq n-2$, and
$\rho(\ora{\hat x_{n-1}x_0})$ for $i=n-1$, respectively.
Then the functional $L$ is just given by $(l_{\eta_1},\dots,l_{\eta_n})$ 
for $\eta_i=\omega_i\eta$. 
After a small variation inside the parallel set $P_\gamma$ we may assume
that the segments $x_0x_k$ and $x_0\hat x_{k+1}$ are regular.
Let $\omega',\omega''\in W$ be so that $\omega'^{-1}\Delta$ 
contains the direction $\rho(\ora{x_0x_k})$ and $\omega''^{-1}\Delta$
contains $\rho(\ora{\hat x_{k+1}x_0})$.

Consider the regular polygon
$q=(x_0,\dots,x_k)\subset A$ and the functional 
$L'=(l_{\eta_0},\dots,l_{\eta_{k-1}},l_{\eta'})$ for
$\eta':=\delta\omega'\eta$. 
That is, $L'$ is the functional 
given in $q$ by taking scalar product with a unit vector in 
the direction $\xi$.
Hence $L'(\sigma(q))=0$.
Set $(\tau_0,\dots,\tau_{k-1},\tau_k):=
(\omega_0,\dots,\omega_{k-1},\delta\omega')$.
Suppose that there are $0\leq i < j\leq k$ such that
$T_{\tau_i}\cap T_{\tau_j}\not\subset \Lambda_\eta^+$.
Corollary~\ref{cor:polygoninapt} and its proof imply that 
there is a polygon $q'=(z_0,\dots,z_k)$ 
with $L'(\sigma(q'))>0$ and with refined side lengths
as near as we want to those of $q$ modulo displacement along $\gamma$.
We can then choose $x_0',x_k'\in P_{\xi,\xi'}$ near $x_0,x_k$ such that 
$x_0'x_k'$ has the same refined side length
(again modulo displacement along $\gamma$) as $z_0z_k$.
The functional $(-l_{\eta'},l_{\eta_{k}},\dots,l_{\eta_{n-1}})$ 
applied to the polygon
$(x_0',x_k',x_{k+1},\dots,x_{n-1})$ is 0 because it is contained in the
parallel set $P_{\xi,\xi'}$.
After displacing the polygon $(x_0',x_k',x_{k+1},\dots,x_{n-1})$ along $\gamma$
we can glue it together to $q'$ and obtain a polygon $p_1$
with $\Delta$-valued side lengths as near as we want to those of $p$
and with $L(\sigma(p_1))>0$
(compare with the proof of Proposition~\ref{prop:unionofcones}). 
This contradicts the assumption at the beginning
of the proof. 
Thus, $T_{\tau_i}\cap T_{\tau_j}\subset \Lambda_\eta^+$
for all $0\leq i < j\leq k$, i.e.\ $(\tau_0.\dots,\tau_k)$
satisfies the property ($\ast$). 
This can be rewritten as 
$T_{\omega_i}\cap T_{\omega_j}\subset \Lambda_\eta^+$ for all
$0\leq i < j \leq k-1$ and
$\bigcup\limits_{i=0}^{k-1}T_{\omega_i} \setminus \Lambda_\eta^+ \subset T_{\omega'}$.

Analogously, considering the polygon $(x_0, x_{k+1},\dots, x_{n-1})$
which is also contained in an apartment in $P_{\xi,\xi'}$
we obtain $T_{\omega_i}\cap T_{\omega_j}\subset \Lambda_\eta^+$ for all
$k+1\leq i < j \leq n-1$
and
$\bigcup\limits_{i=k+1}^{n-1}T_{\omega_i} \setminus \Lambda_\eta^+ \subset T_{\omega''}$.

Consider now the triangle $t=(x_0,x_k,x_{k+1})$
with the functional $L''=(l_{\omega'\eta},l_{\eta_{k}},l_{\omega''\eta})$.
Let $\omega_{k}'\in W$
be so that 
$\omega_{k}'^{-1}\eucch$ contains the direction $\rho(\ora{y\hat x_{k+1}})$.
Then $\eta_{k}=\omega_{k}\eta=\omega_{k}'\eta$.
We want to show that $(\omega',\omega_{k},\omega'')$ 
or $(\omega',\omega_{k}',\omega'')$ have the property ($\ast$).
By Lemma~\ref{lem:crossablewalls} applied to the side $x_0x_k$ we get 
$T_{\omega'}\cap T_{\omega''}, T_{\omega'}\cap T_{\omega_k} \subset \Lambda_\eta^+$.
Again by Lemma~\ref{lem:crossablewalls} now applied to the side 
$x_{k+1}x_0$ we obtain
$T_{\omega''}\cap T_{\omega_k'}\subset \Lambda_\eta^+$.
Therefore if $T_{\omega'}$ or $T_{\omega''}$ is contained in
$\Lambda_\eta^+$, then we are done,
so suppose this is not the case for both of them.
Now by Lemma~\ref{lem:crossablewallsrank2} one of 
$\omega'^{-1}\Delta$, $\omega''^{-1}\Delta$ or $\omega_{k}^{-1}\Delta$ 
must be adjacent to $\rho(\gamma)$.
Notice that for $\omega\in W$, $\omega\Delta$ is adjacent to $\rho(\gamma)$
if and only if $T_\omega\subset \Lambda_\eta^+$ or 
$\Lambda^+\setminus\Lambda_\eta^+\subset T_\omega$.
Since $T_{\omega'}$ and $T_{\omega''}$ are not contained in $\Lambda_\eta^+$
and $T_{\omega'}\cap T_{\omega''}\subset \Lambda_\eta^+$,
then neither $\omega'^{-1}\Delta$ nor $\omega''^{-1}\Delta$ can be
adjacent to $\rho(\gamma)$.
Hence, $\omega_{k}^{-1}\eucch$ must be adjacent to $\rho(\gamma)$.
$T_{\omega'}\cap  T_{\omega_k} \subset \Lambda_\eta^+$ implies that
$T_{\omega_{k}} \subset \Lambda_\eta^+$ and we are also done in this case.

Putting our three considerations above together
we can conclude that $\bar\omega=(\omega_0,\dots,\omega_k,\dots,\linebreak\omega_{n-1})$ 
or $\bar\omega'=(\omega_0,\dots,\omega_k',\dots,\omega_{n-1})$ has the property
($\ast$) and since $p$ is a regular polygon with $L(\sigma(p))=0$, it follows
from Lemma~\ref{lem:noirredundant} that 
$L=L_{\bar\omega\eta}=L_{\bar\omega'\eta}\in {\cal B}_n$.
\qed

Now we are ready to prove our main theorem.

\begin{thm}\label{thm:mainthm}
 Let $X$ be a thick Euclidean building of rank 2. Then $\polspace$ is a convex polyhedral
cone determined by the inequalities $\{L\leq 0\}$ for $L\in {\cal B}_n$. That is,
$$
\polspace = \bigcap_{L\in {\cal B}_n} \{L\leq 0\}.
$$

This inequalities constitute an irredundant set of inequalities.
\end{thm}
\proof
Let $Q\subset {\cal C}_n$ be the subset of open cones such that
$\bigcap\limits_{L\in {\cal B}_n} \{L\leq 0\} = \bigcup\limits_{C\in Q} \bar C$.
Analogously, let $Q'\subset {\cal C}_n$ be the subset of open cones such that
$\polspace = \bigcup\limits_{C\in Q'} \bar C$ (this can be done by
Proposition~\ref{prop:unionofcones}). We have shown in
Lemma~\ref{lem:ineq=weakineq} and Theorem~\ref{thm:weakineq} that $Q'\subset Q$.
Let $C_0\in Q'$ and $C\in Q$. Take a chain $C_0,C_1,\dots,C_k=C \in Q$
such that $\overline{C_i}\cap \overline{C_{i+1}}$ is a face of codimension one.
We prove now inductively that $C_i\in Q'$. Suppose then that
$C_i\in Q'$ and take a regular polygon $p$ with $\sigma(p)$
in the interior of the face $\overline{C_i}\cap \overline{C_{i+1}}$. 
Since $\overline{C_i}\cap \overline{C_{i+1}}$ is not in
the boundary of $\bigcap_{L\in {\cal B}_n} \{L\leq 0\}$, 
it lies in a wall $H_L$ with
neither $L$ nor $-L$ in ${\cal B}_n$. It follows from Lemma~\ref{lem:notpassablewalls}
that $\polspace\cap C_{i+1}$ is not empty and therefore
 $C_{i+1}\subset\polspace$.
Thus $C\in Q'$, and $Q = Q'$.

For $L\in {\cal B}_n$ it is clear that we can find a regular polygon $p$ 
in an apartment $A$
and $\gamma\subset A$ a maximal singular line, 
such that the functional $L$ in $p$
is given by taking scalar product with the direction of $\eta=\gamma(\infty)$. 
In particular, $L(\sigma(p))=0$. 
It is also clear that we can find 
a regular polygon $p'$ in $P_\gamma$ but not contained in any apartment
and such that the functional $L$ in $p'$
is also given by taking scalar product with the direction of $\eta$.
It follows from Lemmata~\ref{lem:wallbyaccident} and \ref{lem:localparallelset}
that $L$ is the only functional in ${\cal B}_n$
 for which it can hold $L(\sigma(p'))=0$.
Thus the inequalities $\{L\leq 0\}$ with $L\in {\cal B}_n$ are irredundant.
\qed

\bigskip\bigskip
\textbf{{\em Acknowledgments.} }
I would like to thank Bernhard Leeb for bringing this problem to my
attention and sharing his ideas in the case of symmetric spaces with me.
A last version of this paper was completed during a stay at the Max-Planck
Institute in Bonn. The author is grateful to the MPI for its financial
support and great hospitality.

\end{document}